\documentclass[11pt]{article}

\usepackage{tikz}
\usepackage{subfigure}
\usepackage{scalefnt}
\usepackage{bm}

\usepackage{amsmath,amsfonts,amssymb,amsthm,enumerate,graphicx}

\usepackage{hyperref}

\usepackage{tikz}
\usetikzlibrary{arrows,backgrounds}
\usetikzlibrary{snakes}
\usetikzlibrary{automata}
\usetikzlibrary{er}
\usetikzlibrary{folding}
\usetikzlibrary{matrix}
\usetikzlibrary{mindmap}
\usetikzlibrary{petri}
\usetikzlibrary{plothandlers}
\usetikzlibrary{plotmarks}
\usetikzlibrary{shapes}
\usetikzlibrary{shapes.symbols}
\usetikzlibrary{shapes.multipart}
\usetikzlibrary{shapes.misc}
\usetikzlibrary{topaths}
\usetikzlibrary{trees}

\usepackage{mathrsfs}
\usepackage{color}
\usepackage{graphics}
\usepackage{epsfig}
\usepackage{tikz,authblk}

\numberwithin{equation}{section}

\newtheorem{theo}{{\bf Theorem}}[section]
\newtheorem{cor}[theo]{{\bf Corollary}}

\newtheorem{defn}[theo]{{\bf Definition}}
\newtheorem{eg}[theo]{{\bf Example}}

\newtheorem{lemma}[theo]{{\bf Lemma}}
\newtheorem{nota}[theo]{{\em Notation}}

\newcommand{\La}{\Lambda}

\newcommand{\si}{\sigma}

\begin{document}


\title{On the structure of cellular pseudomanifolds}

\author[1] {Bhaskar Bagchi\thanks{A retired professor of Indian Statistical Institute,  Bangalore, India.
}}
\author[2] {Basudeb Datta}

\affil[1]{B36, Good Earth, Malhar Medley, Bangalore 560\,074, India. \newline Email: {\it bhaskarbagchi53@gmail.com}.}

\affil[2]{Department of Mathematics, Indian Institute of Science, Bangalore 560\,012, India. \newline Institute for Advancing Intelligence, TCG CREST,  Kolkata 700\,091, India. \newline 
Emails: {\it dattab@iisc.ac.in,  bdatta17@gmail.com}.}

\date{July 04, 2023}

\maketitle

%

\vspace{-5mm}




\bigskip

\hrule

\begin{abstract}
In this paper we study the structure of cellular pseudomanifolds (aka abstract polytopes). These are natural combinatorial generalisations of polytopal spheres (i.e., boundary complexes of convex polytopes). This class is closed under natural notions of duality and product. We show that they are also closed under an operation of direct product. Any cellular pseudomanifold and it's dual have homeomorphic geometric carriers, while the geometric carrier of the product of two of them is homeomorphic to the product of the carriers of the factors. The excess of a cellular pseudomanifold is defined as the non-negative integer $n-d-2$ where $d$ is the dimension and $n$ is the number of vertices. We completely classify the cellular pseudo manifolds of excess $\leq 1$, and make some progress towards classifying those of excess 2.
\end{abstract}

\tableofcontents

\newpage

\section{\sc Introduction}

We define a (simplicial) {\em normal pseudomanifold} to be a finite collection ${\cal F}$ of finite sets satisfying the axioms $(0), \dots, (4)$ below. The set $V = V({\cal F}) := \cup_{F\in {\cal F}} F$ is called the {\em vertex-set} of ${\cal F}$ and the elements of $V$ are called the {\em vertices} of ${\cal F}$. The elements of ${\cal F}$ are called the {\em faces} of ${\cal F}$.
The faces other than $\emptyset$ and $V$ are called the {\em proper faces}. If a proper face $F$ of ${\cal F}$ is not properly contained in any proper face of $\cal F$ then $F$ is called a {\em facet} of $\cal F$. $\La(\cal F)$ denote the graph whose vertices are the facets of $\cal F$, and two facets $F_1\neq F_2$ are adjacent in $\La(\cal F)$ if $\#(F_1 \setminus F_2) = 1 = \#(F_2 \setminus F_1)$. The axioms satisfied by $\cal F$ are as follows\,:
\begin{enumerate}[(1)] 
\item[$(0)$] 
({\em Trivial faces}) $\emptyset\in {\cal F}$, $V\in {\cal F}$ and $\{v\}\in {\cal F}$ for all $v\in V$.
\item[$(1)$] 
({\em Simplicial}) Every subset of a proper face   of $\cal F$ is again a face of $\cal F$.
\item[$(2)$] 
({\em Dimension}) There is an integer $d\geq -1$    such that $\#(F) = d+1$ for every facet of $\cal F$. This number is called the {\em dimension} of $\cal F$ and is denoted by $\dim(\cal F)$. 
\item[$(3)$] 
({\em Strong connectedness}) The graph $\La(\cal F)$ is   connected.
\item[$(4)$] 
({\em Induction}) If $\dim({\cal F}) \leq  0$ then   $\#(V({\cal F})) = \dim({\cal F}) + 2$. If $\dim({\cal F}) > 0$  then for each vertex $x\in V(\cal F)$, ${\rm lk}_{\cal F} (x) := \{F     \setminus \{x\} ~ \colon ~ x\in F\in {\cal F}\}$ is a normal pseudomanifold.
\end{enumerate}

This definition of a normal pseudomanifold differs from the usual definitions found in the literature. The difference is that  we have insisted on taking the entire vertex-set $V$ as a face. Actually, we regard $V$ as a kind of ``fake" face which affects nothing; in particular the
definition of the geometric carrier given below is un-affected by this change. However, the introduction of $V$ as a face in $\cal F$ converts the poset $({\cal F}, \subseteq)$ into a (complete) lattice and this simplifies many of the definitions and proofs which occur in this paper. Indeed, the normal pseudomanifold $\cal F$ is completely determined by the associated poset $({\cal F}, \subseteq)$, prompting us to define the main objects of study in
this paper - certain generalizations of normal pseudomanifolds called cellular manifolds - as a class of lattices. The notion of a (simplicial) pseudomanifold  is obtained from the above definition by omitting the last part of Axiom (4).  

Two normal pseudomanifolds $\cal F$ and $\cal F^{\,\prime}$ are said to be {\em isomorphic} if there is a bijection $f\colon {\cal F} \to {\cal F^{\,\prime}}$ which preserves inclusion. Notice that such an isomorphism is determined by its restriction to $V(\cal F)$. (The latter is a bijection from $V(\cal F)$ onto $V(\cal F^{\,\prime})$.) We identify two normal pseudomanifolds if they are isomorphic.

The simplest example of a normal pseudomanifold is the standard $d$-sphere ($d\geq -1$) $S^{\,d}_{d+2}$ whose faces are all the subsets of a set of size $d+2$. Notice that for any proper face $B$ of a normal pseudomanifold $\cal F$, the interval $[\emptyset, B] := \{F \in
{\cal F} ~ \colon ~ F\subseteq B\}$ is a standard sphere. The normal pseudomanifold $[\emptyset, B]$ is called the {\em boundary} of $B$ and is denoted by $\partial B$. On the other hand, for any proper face $A$ of $\cal F$, the interval $[A, V] := \{F \in {\cal F} ~ \colon ~ A\subseteq F\}$ is isomorphic as a poset to  ${\rm lk}_{\cal F}(A) := \{F \setminus A ~ \colon ~
A \subseteq F \in {\cal F}\}$. The normal pseudomanifold ${\rm lk}_{\cal F}(A)$ is called the {\em link} of $A$ in $\cal F$. It follows that, more generally, whenever $A \subseteq B$ are two faces of a normal pseudomanifold $\cal F$, the interval $[A, B] := \{F \in {\cal F} ~ \colon A \subseteq F \subseteq B\}$, equipped with the partial order $\subseteq$, is the poset of a normal pseudomanifold. This observation forms the motivation behind the definition of a
cellular pseudomanifold. To gain an informal understanding of this generalization, notice that
$S^{\,-1}_1$ and $S^{\,0}_2$ are the only  normal pseudomanifolds of dimensions $\leq 0$ and, for integers $n\geq 3$, the $n$-cycles $S^{\,1}_n$ (i.e., the unique $n$-vertex connected graph which is regular of degree 2) are the only 1-dimensional  normal pseudomanifolds. Indeed, our generalization is such that these are also the only cellular pseudomanifolds of dimension $\leq 1$. However, by definition, a  normal pseudomanifold of dimension 2 has only triangles ($S^{\,1}_3$'s) as facets while the general 2-dimensional cellular pseudomanifold permits arbitrary polygons ($S^{\,1}_n$ with possibly varying gonality $n$) as facets. For instance, the boundary complex of the icosahedron is a  normal pseudomanifold while the boundary complex of the dodecahedron is a cellular pseudomanifold which is not simplicial (its facets are pentagons).
As indicated above (and elaborated in the next section), the boundary complex of a convex
polytope is a special example of cellular pseudomanifolds. In fact, this paper is partly motivated by the desire to classify the $d$-polytopes on few vertices up to cellular isomorphism.

Recall that a poset (partially ordered set) $(X, <)$ is called a
{\em lattice} if for any two elements $x$, $y$ in $X$, there is a
(unique) meet $x\wedge y\in X$ and a (unique) join $x\vee y\in X$
such that $x\wedge y$ is the largest element smaller than both $x$
and $y$, while $x\vee y$ is the smallest element larger than both
$x$ and $y$. In particular, a finite lattice $(X, <)$ has a unique
smallest element, usually denoted by {\bf 0}, and a unique largest
element, usually denoted by {\bf 1}. For two elements $a$, $b$ in
a lattice $(X, <)$, we write $a\leq b$ if $a = b$ or $a < b$. Two
elements $a$, $b$ are said to be {\em incident} if $a = b$, $a <
b$ or $b < a$. If $a < b$ are two elements in a lattice $(X, <)$
then the {\em interval} $[a, b] := \{x \in X ~ : ~ a \leq x\leq
b\}$ is clearly a lattice with the induced partial order. A subset
$C$ of $X$ is called a {\em chain} if any two elements in $C$ are
incident.

A finite lattice $(X, <)$ is said to be {\em ranked} if there is a
rank function $\rho \, \colon X\to \mathbb{N}$ such that for each $x\in
X$, all the maximal chains ${\bf 0} = x_0 < \cdots < x_r = x$
joining {\bf 0} to $x$ have length $r=\rho(x)$. Thus $\rho({\bf
0}) = 0$. For a ranked lattice $L=(X, <)$, let $\Lambda(L)$ denote
the graph whose vertices are elements in $X$ of rank $\rho({\bf
1}) -1$, where two elements $\sigma$ and $\gamma$ (of ranks
$\rho({\bf 1}) - 1$) are adjacent if and only if $\rho(\sigma
\wedge \gamma) = \rho({\bf 1}) - 2$.

\begin{defn}  {\rm A ranked lattice $M = (X,
<)$ is said to be a {\em cellular pseudomanifold} if for any two
elements $x$, $z$ in $X$ with $x < z$ (i) there are exactly two
elements $y$ satisfying $x < y < z$ wherever $\rho(z) = \rho(x) +
2$ and (ii) $\Lambda([x, z])$ is connected whenever $\rho(z) >
\rho(x) +2$. }
\end{defn}

Cellular pseudomanifolds are also known as {\em abstract polytopes} in the literature (cf. \cite{ms}), but we prefer the first name since cellular pseudomanifolds which are isomorphic to the boundary complexes of convex polytopes form a proper subclass of the class of cellular pseudomanifolds.


The class of cellular pseudomanifolds has certain conceptual advantages over the class of 
normal pseudomanifolds. Unlike the latter, the former class is closed under a natural notion of duality as well as under a natural product. A cellular pseudomanifold and its dual have homeomorphic geometric carriers. The geometric carrier of the product of two cellular manifolds is the topological product of the geometric carriers of the factors. The notion of join in the category of normal pseudomanifolds has a natural extension to cellular
pseudomanifolds - and, in fact, the extension looks even simpler than the original definition.

The elements of a cellular pseudomanifold $M$ are called the {\em faces} of $M$. For any face $f$, $\dim(f) := \rho(f) -1$ is called the {\em dimension} of $f$. The integer $\dim(M) := \dim({\bf 1}) - 1$ is called the {\em dimension} of $M$. The faces of dimensions 0 and 1 are
called the {\em vertices} and {\em edges} respectively. A face of dimension $i$ is said to be an
{\em $i$-face}. The $\dim(M)$-faces are also called the facets. The faces ${\bf 0}$ and ${\bf
1}$ of a cellular pseudomanifold $M$ are also called the {\em trivial faces} of $M$. Thus, $S^{\,-1}_1$ is the cellular pseudomanifold without proper faces. 

If $M_1 = (X_{1},<_1)$ and $M_2 = (X_{2},<_2)$ are two cellular pseudomanifolds, then a {\em cellular isomorphism} from $M_1$ to $M_2$ is a bijection $\pi : X_{1} \to X_{2}$ such that $f<_1 g$ if and only if $\pi(f)<_2 \pi(g)$. The cellular pseudomanifolds $M_{1}, M_{2}$ are called {\em isomorphic} when such an isomorphism exists. We identify two cellular
pseudomanifolds if they are isomorphic.

\section{\sc Cellular pseudomanifolds}

For a normal pseudomanifold $\cal F$, if $V({\cal F}) = \{v_1, \dots,
v_n\}$ then consider an $(n-1)$-simplex with vertices $\{x_1,
\dots, x_n\}$ in $\mathbb{R}^{n-1}$. If $F = \{v_{i_1}, \dots, v_{i_k}\}$
is a proper face of $\cal F$ then the convex set spanned by
$\{x_{i_1}, \dots, x_{i_k}\}$ is said to be the {\em geometric
carrier} of $F$ and is denoted by $|F|$. The space ($\subseteq
\mathbb{R}^{n-1}$) obtained by taking the union of the geometric carriers
of proper faces of $\cal F$ is called the {\em geometric carrier}
of $\cal F$ and is denoted by $|{\cal F}|$.

By a {\em subdivision} of a normal pseudomanifold $\cal F$ we mean a
normal pseudomanifold $\cal F^{\,\prime}$ together with a homeomorphism
$f \colon |{\cal F}^{\,\prime}| \to |{\cal F}|$ such that the
restriction of $f$ on $|A|$ is linear for each proper face $A$ of
$\cal F^{\,\prime}$. Two normal pseudomanifolds ${\cal F}_1$ and ${\cal
F}_2$ are called {\em combinatorially equivalent} (denoted by
${\cal F}_1 \approx {\cal F}_2$) if they have isomorphic
subdivisions. So, ${\cal F}_1 \approx {\cal F}_2$ if and only if
$|{\cal F}_1|$ and $|{\cal F}_2|$ are pl homeomorphic (\cite{rs}).
A normal pseudomanifold $\cal F$ is called a {\em combinatorial
$d$-sphere} if it is combinatorially equivalent to $S^{\,d}_{d +
2}$. (Thus the standard spheres and $n$-cycles are the simplest
examples of combinatorial spheres.) A $d$-dimensional
normal pseudomanifold is called a {\em combinatorial $d$-manifold} if the
link of each vertex in it is a combinatorial $(d - 1)$-sphere.

For any cellular pseudomanifold $M$, $V(M)$ will denote the set of
vertices of $M$. For any face $f$ of $M$, let ${\rm sh}(f)$, the
{\em shadow} of $f$, denote the set of all vertices $v \in V(M)$
such that $v \leq f$. It is easy to see by induction on $\dim(f)$
that $f$ is uniquely determined by ${\rm sh}(f)$; indeed, $f$ is
the join of all the vertices in ${\rm sh}(f)$. In consequence, for
any two faces $f$ and $g$, $f < g$ if and only if ${\rm sh}(f)
\subset {\rm sh}(g)$. {\em Thus we may identify each face with its
shadow, hence identifying $M$ with a family of subsets of $V(M)$,
ordered by set inclusion \/}. We shall use this identification in
many of the proofs in this paper. From this identification, it is
immediate that the shadow of any face $f$ contains at least
$\dim(f) + 1$ vertices. From our identification it is also clear
that equality holds here (for all faces) if and only if $M$ is a normal 
pseudomanifold.

For two faces $f < g$ of a cellular pseudomanifold $M$, the {\em
interval} $[f, g]$ is clearly a cellular pseudomanifold (with the
induced partial order) of dimension $= \dim(g) - \dim(f) - 2$. The
interval $[{\bf 0}, g]$ is called the {\em boundary} of $g$,
sometimes denoted by $\partial_M g$ (or $\partial g$). The
interval $[f, {\bf 1}]$ is called the {\em link} (in $M$) of $f$,
and is denoted by ${\rm lk}_M(f)$ (or ${\rm lk}(f)$). Note that,
in particular, $\partial_M {\bf 1} = M = {\rm lk}_M({\bf 0})$.

 An isomorphism from a
cellular pseudomanifold $M$ to itself is called an {\em
automorphism} of $M$. Note that any automorphism of $M$ restricts
to a permutation of $V(M)$, which in turn determines the
automorphism. The automorphisms of $M$ form a group under
composition. It is called the {\em automorphism group} of $M$, and
is denoted by ${\rm Aut}(M)$.

For a $d$-dimensional cellular pseudomanifold $M$, length of any
maximal chain from ${\bf 0}$ to ${\bf 1}$ is $d+2$. This implies
that each face of $M$ is incident to a facet of $M$. Here we
have\,:

\begin{lemma}  \label{l2.1} Let $M$ be a cellular
pseudomanifold and let $\alpha$ be a face of $M$. If a vertex $v$
is not incident with $\alpha$ then there exists a facet $\sigma$
incident with $\alpha$ but not with $v$.
\end{lemma}

\noindent {\em Proof.} Backward induction on $\dim(\alpha)$. Since
$v$ is not incident with $\alpha$, $\alpha\neq {\bf 1}$. The
result is obvious for $\dim(\alpha) = \dim(M)$. So, assume
$\dim(\alpha) < \dim(M)$ and $v$ is a vertex not incident with
$\alpha$. Let $\beta_1$ and $\beta_2$ be two faces incident with
$\alpha$ such that $\dim(\beta_1) = \dim(\beta_2) = \dim(\alpha) +
1$. If $v$ is incident with both $\beta_1$, $\beta_2$ then $\alpha
< \beta_1 \wedge \beta_2 < \beta_1$. This is not possible since
$\dim(\beta_1) = \dim(\alpha) + 1$. So, $v$ is not incident with
at least one of $\beta_1$, $\beta_2$, say $v$ is not incident with
$\beta_1$. Then by induction hypothesis, there exists a facet
$\sigma$ incident with $\beta_1$ but not with $v$. Then $\sigma$
is incident with $\alpha$ but not with $v$. This completes the
proof. \hfill $\Box$

For any cellular pseudomanifold $M$, let $\Lambda^{\ast}(M)$
denote the graph whose vertices are vertices of $M$, where $u,
v\in V(M)$ are adjacent if $\dim(u\vee v) = 1$. Observe that the edge graph of $M$ (graph consists of vertices and edges of $M$) is isomorphic to $\Lambda^{\ast}(M)$ via the map $e\mapsto {\rm sh}(e)$ for any edge $e$ of $M$. We have\,:

\begin{lemma}  \label{l2.2} For any cellular
pseudomanifold $M$ of dimension $\geq 1$, $\Lambda^{\ast}(M)$ is
connected. In consequence, whenever $x< y$ are faces with $\dim(y)
> \dim(x) +2$, $\Lambda^{\ast}([x, y])$ is connected.
\end{lemma}

\noindent {\em Proof.} The proof is by induction on $\dim(M)$. If
$\dim(M)=1$ then $\Lambda(M)$ is a cycle and hence
$\Lambda^{\ast}(M)$ is a cycle. Thus, the result is true for
dimension $= 1$. So, let $\dim(M)> 1$ and suppose that we have the
result for all cellular pseudomanifolds of dimension $< \dim(M)$.

Since $\Lambda(M)$ is connected, it is sufficient to show that for
any two facets $\alpha$, $\beta$ which are adjacent in
$\Lambda(M)$ and for any two vertices $u$, $v$ with $u<\alpha$,
$v< \beta$ we have a path from $u$ to $v$ in $\Lambda^{\ast}(M)$.
Let $w$ be a vertex incident with $\alpha\wedge\beta$. Then $u$,
$w$ are vertices in $\partial \alpha$ and $\dim(\partial \alpha) <
\dim(\alpha)=\dim(M)$. Therefore, by induction hypothesis, there
is a path in $\Lambda^{\ast}(M)$ joining $u$ to $w$. Similarly,
there is a path in $\Lambda^{\ast}(M)$ joining $w$ to $v$.
Juxtaposing these two paths, we obtain a path from $u$ to $v$ in
$\Lambda^{\ast}(M)$. This completes the proof. \hfill $\Box$

\begin{defn}  {\rm For any cellular pseudomanifold $M
= (X, <)$, the opposite poset $M^{\ast} = (X, <^{\ast})$ is called
the {\em dual} of $M$. Thus, the faces of $M^{\ast}$ are precisely
the faces of $M$. For faces $x$ and $y$, $x <^{\ast} y$ in
$M^{\ast}$ if and only if  $y < x$ in $M$. }
\end{defn}

From Lemma \ref{l2.2}, we have the following.

\begin{cor}  \label{c2.1} For a cellular pseudomanifold
$M$, the dual $M^{\ast}$ is a cellular pseudomanifold.
\end{cor}

\begin{eg}  \label{polytope} {\rm Recall that a set
$X \subseteq \mathbb{R}^{d}$ is called a {\em polytope} of dimension $d$
($\geq 0$) if $X$ is the convex hull of finitely many points and
$X$ has a non-empty interior. A closed face in $X$ is the
intersection of $X$ with a hyperplane $H$ such that $X$ is
contained in one of the closed half-spaces determined by $H$. We
also admit $X$ itself as a closed face of $X$. Then the poset
$\partial X$ consisting of all the closed faces of $X$, ordered by
set inclusion, is an example of a $(d-1)$-dimensional cellular
pseudomanifold. It is called the {\em boundary complex} of $X$. If
the number of vertices of $X$ is $d + 1$ (equivalently, $X$ is a
$d$-simplex) then $\partial X$ is normal pseudomanifold and is
called the {\em standard $(d-1)$-sphere}. The standard
$(d-1)$-sphere on the vertex set $V$ is denoted by $S^{\,d-1}_{d +
1}(V)$ (or simply by $S^{\,d-1}_{d + 1}$). If a cellular
pseudomanifold is isomorphic to the boundary complex of a polytope
then we say that it is {\em polytopal}. }
\end{eg}

For a cellular pseudomanifold $M$, let $B(M)$ be the set of all the chains 
in $M\setminus\{{\bf 0}, {\bf 1}\}$ together with $M
\setminus \{{\bf 0}, {\bf 1}\}$. Then $B(M)$ is a normal 
pseudomanifold and is called the {\em barycentric subdivision} of
$M$. (So, the vertices of $B(M)$ are the proper faces of $M$.) A
$d$-dimensional cellular pseudomanifold $M$ is called a {\em
cellular $d$-sphere} (respectively, {\em cellular $d$-manifold})
if $B(M)$ is a combinatorial $d$-sphere (respectively,
combinatorial $d$-manifold). So, a polytopal $d$-sphere is a
cellular $d$-sphere.

The geometric carrier of $B(M)$ is called the {\em geometric
carrier} of $M$ and is denoted by $|M|$. For a proper face $f$ in
$M$, let $B(f) := \{\sigma\cap [{\bf 0}, f] ~ \colon ~ \sigma\in
B(M)\}\subset B(M)$ and let $|f| := \cup_{\sigma\in B(f)} |\sigma|
\subseteq |B(M)|$. Then $|f|$ is called the {\em geometric
carrier} of $f$. Clearly, $|M|$ is the union of the geometric
carriers of all its facets. Thus, $|S^{\,-1}_1| = \emptyset$ and
for $d\geq 0$ the geometric carrier of a cellular $d$-sphere is
homeomorphic to the sphere $S^{\,d}$.

Recall that a set $P \subseteq \mathbb{R}^n$ is a {\em compact
polyhedron} if $P$ can be express as a finite union of geometric
simplices (cf. \cite{rs}). So, a polytope is a compact polyhedron.
More generally, if $M$ is a $d$-dimensional cellular
pseudomanifold and $B(M)$ has a geometric carrier in $\mathbb{R}^n$ (for
some $n$) then $|M|$ is a compact polyhedron in $\mathbb{R}^n$.

\begin{lemma}  \label{l2.3}
A cellular pseudomanifold $M$ is a cellular manifold if and only
if $[\alpha, \beta]$ is a cellular sphere for every interval
$[\alpha, \beta]\neq M$.
\end{lemma}

\noindent {\em Proof.} Let $\sigma$ is a proper face of $M$. Then
$\sigma$ is a vertex of $B(M)$ and any proper face $A$ containing
$\sigma$ in $B(M)$ is of the form $\{\sigma_0, \dots, \sigma_{i},
\dots, \sigma_{i+j}\}$ for some faces $\sigma_0, \dots,
\sigma_{i+j}$ of $M$ with $\sigma_0 < \cdots < \sigma_{i-1} <
\sigma = \sigma_{i} < \cdots < \sigma_{i+j}$. Then $A\setminus\{\sigma\} =
\{\sigma_0, \dots, \sigma_{i-1}\}\cup \{\sigma_{i+1}, \dots, \sigma_{i
+ j}\} \in B(\partial_M \sigma) \ast B({\rm lk}_{M}(\sigma))$.
Conversely, any face of $B(\partial_M \sigma) \ast B({\rm
lk}_{M}(\sigma))$ is in ${\rm lk}_{B(M)}(\sigma)$. Thus
   \begin{eqnarray} \label{link}
   {\rm lk}_{B(M)}(\sigma) & = & B(\partial_M \sigma) \ast B({\rm
   lk}_{M}(\sigma)).
\end{eqnarray}

Let $M$ be a cellular manifold, $\sigma$ be a proper face of $M$
and $L:= {\rm lk}_{B(M)}(\sigma)$. Then $B(M)$ is a combinatorial
manifold. So, $L$ is a combinatorial sphere (and hence a
combinatorial manifold). Thus, if $\beta$ is a facet of
$B(\partial_M \sigma)$ then, by (\ref{link}), ${\rm lk}_{L}(\beta)
= B({\rm lk}_{M}(\sigma))$. This implies that $B({\rm
lk}_{M}(\sigma))$ is a combinatorial sphere. Similarly,
$B(\partial_M \sigma)$ is a combinatorial sphere. Thus, ${\rm
lk}_M(\sigma)$ ($= [\sigma, {\bf 1}]$) and $\partial_M(\sigma)$
($= [{\bf 0}, \sigma]$) are cellular spheres.

Now, let ${\bf 0} < \alpha < \beta < {\bf 1}$. Let $N : = [{\bf
0}, \beta]$. Then, by above, $N$ is a cellular sphere (and hence a
cellular manifold). Since $[\alpha, \beta] = {\rm lk}_N(\alpha)$,
by above, $[\alpha, \beta]$ is a cellular sphere. Thus, if $M$ is
a cellular manifold then $[\alpha, \beta]$ is a cellular sphere
for each interval $[\alpha, \beta]\neq M$.

Conversely, suppose $[\alpha, \beta]$ is a cellular sphere for
each interval $[\alpha, \beta]\neq M$. Then, $[{\bf 0}, \sigma]$
and $[\sigma, {\bf 1}]$ cellular spheres for any proper face of
$M$. Therefore, by (\ref{link}), ${\rm lk}_{B(M)}(\sigma)$ is a
combinatorial sphere for each vertex $\sigma$ of $B(M)$. This
implies that $B(M)$ is a combinatorial manifold. This proves the
lemma. \hfill $\Box$

\medskip

It is trivial that for the existence of a $d$-dimensional cellular
pseudomanifold on $n$ vertices, we must have $n \geq d + 2$.
Equality holds here only for the standard $d$-sphere
$S^{\,d}_{d+2}$. This motivates the following.

\begin{defn} 
{\rm The {\em excess} $e(M)$ of a $d$-dimensional cellular
pseudomanifold $M$ on $n$ vertices is the non-negative number
$e(M) = n - d - 2$.}
\end{defn}

\begin{defn} 
{\rm Let $M_1=(X_1, <_1)$ and $M_2 = (X_2, <_2)$ be two cellular
pseudomanifolds. 
\begin{enumerate}[(i)]
\item
The {\em direct product} $M_1 \otimes M_2 = (X_1
\times X_2, <)$ is the poset defined by $(x_1, x_2) < (y_1, y_2)$
if and only if either $x_1 \leq_1 y_1$ and $x_2 <_2 y_2$ or $x_1
<_1 y_1$ and $x_2 \leq_2 y_2$. 
\item 
The {\em join} $M_1 \ast M_2$ is
the sub-poset of $M_1 \otimes M_2$ consisting of all the faces
$(\alpha, \beta)$ where either both of $\alpha$, $\beta$ are $<
{\bf 1}$ or else both are $ = {\bf 1}$. 
\item 
The {\em cartesian
product} $M_1 \times M_2$ is the sub-poset of $M_1 \ast M_2$
consisting of all the faces $(\alpha, \beta)$ where either both of
$\alpha$, $\beta$ are $> {\bf 0}$ or else both are $ = {\bf 0}$.
\item[]
So, $M_1 \ast M_2 = (((X_1 \setminus \{{\bf 1}\}) \times (X_2
\setminus \{{\bf 1}\})) \cup \{({\bf 1}, {\bf 1})\}, <)$ and $M_1
\times M_2 = (((X_1 \setminus \{{\bf 0}, {\bf 1}\}) \times (X_2
\setminus \{{\bf 0}, {\bf 1}\}))\cup \{({\bf 0}, {\bf 0}), ({\bf
1}, {\bf 1})\}, <)$. 
\end{enumerate}
}
\end{defn}

\begin{defn}  \label{star}
{\rm Let $A \subseteq \mathbb{R}^m$ and $B \subseteq \mathbb{R}^n$ be two
compact polyhedra. Let $i \colon \mathbb{R}^{m} \to \mathbb{R}^{m + n + 1}$, $j
\colon \mathbb{R}^{n} \to \mathbb{R}^{m + n + 1}$ be the maps given by $i(x_1,
\dots, x_m) = (x_1, \dots, x_m, 0, \dots, 0)$, $j(y_1, \dots, y_n)
= (0, \dots, 0, y_1, \dots, y_n, 1)$. Let
$$
A \ast B := \{tx + (1-t)y ~ : ~ x \in i(A), ~ y \in j(B), t \in
[0, 1]\}.
$$
Then $A \ast B$ is a polyhedron in $\mathbb{R}^{m+n+1}$ and is called the
(external) {\em join} of $A$ and $B$ (cf. \cite{rs}). }
\end{defn}

If $K$ and $L$ are normal pseudomanifolds with $|K| \subseteq
\mathbb{R}^m$ and $|L| \subseteq \mathbb{R}^n$ then it is easy to see that the
polyhedron $|K| \ast |L| \subseteq \mathbb{R}^{m+n+1}$ is pl homeomorphic
to $|K \ast L|$. If $K$ and $L$ are cellular pseudomanifolds then
also same is true by Theorem \ref{t2.2} $(i)$ below. Here we
have\,:

\begin{lemma}  \label{l2.4} If $X \subseteq
\mathbb{R}^m$ is an $m$-polytope and $Y \subseteq \mathbb{R}^n$ is an
$n$-polytope then $X\ast Y \subseteq \mathbb{R}^{m + n +1}$ is an $(m + n
+ 1)$-polytope with $\partial (X \ast Y) = (\partial X) \otimes
(\partial Y)$. So, the direct product of the boundary complexes of
two polytopes is the boundary complex of their join.
\end{lemma}

\noindent {\em Proof.} Let $i \colon \mathbb{R}^{m} \to \mathbb{R}^{m + n + 1}$,
$j \colon \mathbb{R}^{n} \to \mathbb{R}^{m + n + 1}$ be as in Definition
\ref{star}. For $U \subseteq \mathbb{R}^{m + n +1}$, let $\langle U
\rangle$ denote the convex span of $U$. Then $X \ast Y = \langle
i(X) \cup \, j(Y) \rangle$.

If $a$ is in the interior of $X$ and $b$ is in the interior of $Y$
then the set $\{ta + (1-t)b ~ : ~ 0 < t < 1\}$ is contained in the
interior of $X\ast Y$. Thus, $X \ast Y$ is an $(m + n +
1)$-polytope. Now, if $P$ is a hyperplane in $\mathbb{R}^{m + n + 1}$
which touches $X \ast Y$ then $(P\cap i(\mathbb{R}^m)) \cap (X \ast Y) =
\emptyset$, $i(X)$ or $i(A)$ for some proper closed face $A$ of
$X$. Similarly, $(P \cap j(\mathbb{R}^n)) \cap (X \ast Y) = \emptyset$,
$j(Y)$ or $j(B)$ for some proper closed face $B$ of $Y$. Hence the
boundary complex $\partial (X \ast Y)$ is $\{\emptyset, \, W, \,
i(X), \, j(Y), \, i(A), \, j(B), \, \langle i(A) \cup \, j(B)
\rangle = A\ast B, \, \langle i(X) \cup \, j(B) \rangle = X \ast
B, \, \langle i(A) \cup \, j(Y) \rangle = A \ast Y \, : \,
\emptyset \neq A < X, \, \emptyset \neq B < Y\}$. Therefore,
$\partial (X \ast Y)$ is isomorphic to $(\partial X) \otimes
(\partial Y)$. \hfill $\Box$

\begin{lemma}  \label{l2.5}
If $M_1$ and  $M_2$ are cellular pseudomanifolds of dimensions
$\geq 0$ then $M_1 \otimes M_2$, $M_1 \ast M_2$ and $M_1 \times
M_2$ are cellular pseudomanifolds. Moreover, $\dim(M_1 \otimes
M_2) - 2 = \dim(M_1 \ast M_2) - 1 = \dim(M_1 \times M_2) =
\dim(M_1) + \dim(M_2)$, $\#(V(M_1 \otimes M_2)) = \#(V(M_1 \ast
M_2)) = \#(V(M_1)) + \#(V(M_2))$ and $\#(V(M_1 \times M_2)) =
\#(V(M_1))\times \#(V(M_2))$. In consequence, $e(M_1 \otimes M_2)
= e(M_1 \ast M_2) -1 = e(M_1) + e(M_2)$.
\end{lemma}

\noindent {\em Proof.} First consider the lattice $M_1\otimes
M_2$. Let $(\alpha_1, \beta_1) < (\alpha_2, \beta_2)$ and
$\rho((\alpha_1, \beta_1)) = \rho((\alpha_2, \beta_2)) - 2$. Then
either (i) $\alpha_1=\alpha_2$, $\beta_1 < \beta_2$,
$\rho(\beta_1) = \rho(\beta_2) - 2$, (ii) $\alpha_1 < \alpha_2$,
$\beta_1 = \beta_2$, $\rho(\alpha_1) = \rho(\alpha_2) - 2$ or
(iii) $\alpha_1 < \alpha_2$, $\beta_1 < \beta_2$, $\rho(\alpha_1)
= \rho(\alpha_2) - 1$, $\rho(\beta_1) = \rho(\beta_2) - 1$. In the
first case, if $\gamma_1$, $\gamma_2$ lie between $\beta_1$ and
$\beta_2$ then $(\alpha_1, \gamma_1)$, $(\alpha_1, \gamma_2)$ are
the only faces between $(\alpha_1, \beta_1)$ and $(\alpha_2,
\beta_2)$. The second case is similar to the first case. In the
last case, $(\alpha_1, \beta_2)$ and $(\alpha_2, \beta_1)$ are the
only faces between $(\alpha_1, \beta_1)$ and $(\alpha_2,
\beta_2)$. Thus, in all the cases, there are exactly two faces
between $(\alpha_1, \beta_1)$ and $(\alpha_2, \beta_2)$.

Let $(\alpha_1, \beta_1) < (\alpha_2, \beta_2)$ and
$\rho((\alpha_1, \beta_1)) \leq \rho((\alpha_2, \beta_2)) - 3$.
Clearly, if $(\gamma, \delta)$ is a facet of $[(\alpha_1,
\beta_1), (\alpha_2, \beta_2)]$ then either $\gamma = \alpha_2$
and $\delta$ is a facet of $[\beta_1, \beta_2]$ or $\gamma$ is a
facet of $[\alpha_1, \alpha_2]$ and $\delta = \beta_2$. Thus, to
show that $\La([(\alpha_1, \beta_1), (\alpha_2, \beta_2)])$ is
connected, it is enough to show that (i) there is a path between
$(\alpha_2, \delta_1)$ and $(\alpha_2, \delta_2)$, where
$\delta_1$, $\delta_2$ are facets of $[\beta_1, \beta_2]$ and (ii)
there is a path between $(\alpha_2, \delta)$ and $(\gamma,
\beta_2)$, where $\delta$ is a facet of $[\beta_1, \beta_2]$ and
$\gamma$ is a facet of $[\alpha_1, \alpha_2]$. In the first case,
let $\delta_1 = \eta_1, \dots, \eta_k = \delta_2$ be a path in
$\La([\beta_1, \beta_2])$. Then $(\alpha_2, \delta_1) = (\alpha_2,
\eta_1), \dots, (\alpha_2, \eta_k) = (\alpha_2, \delta_2)$ is a
path in $\La([(\alpha_1, \beta_1), (\alpha_2, \beta_2)])$. In the
last case, $(\alpha_2, \delta) \wedge (\gamma, \beta_2) = (\gamma,
\delta)$ and $\rho((\gamma, \delta)) = \rho((\alpha_2, \delta)) -
1$. Therefore, $(\alpha_2, \delta)$ and $(\gamma, \beta_2)$ are
adjacent in $\La([(\alpha_2, \beta_1), (\alpha_2, \beta_2)])$.
Thus, $M_1 \otimes M_2$ is a cellular pseudomanifold. Same
arguments show that $M_1 \ast M_2$ and $M_1 \times M_2$ are
cellular pseudomanifolds. The rest follow from the definitions.
\hfill $\Box$

\bigskip

From the definitions we get\,:

\begin{lemma}  \label{l2.6}
Let $M$ and  $N$ be two cellular pseudomanifolds of dimensions
$\geq 0$. Let $[\mbox{-} , \mbox{-}]_{\otimes}$, $[\mbox{-} ,
\mbox{-}]_{\ast}$ and $[\mbox{-} , \mbox{-}]_{\times}$ denote the
intervals in $M \otimes N$, $M \ast N$ and $M \times N$
respectively. If $\alpha \in M$, $\beta \in N$ are proper faces
then $[({\bf 0}, {\bf 0}), (\alpha, \beta)]_{\otimes} = [{\bf 0},
\alpha] \otimes [{\bf 0}, \beta] = [({\bf 0}, {\bf 0}), (\alpha,
\beta)]_{\ast}$ and hence $[({\bf 0}, {\bf 0}), (\alpha,
\beta)]_{\times} = ([{\bf 0}, \alpha] \otimes [{\bf 0}, \beta])
\cap (M \times N)$. In consequence, if $f$ is a facet of $\partial
(\alpha, \beta)$ $($in $M \otimes N$, $M \ast N$ or $M \times N)$
then either $f = (g, \beta)$ for some facet $g$ of $\partial
\alpha$ or $f = (\alpha, h)$ for some facet $h$ of $\partial
\beta$.
\end{lemma}

\begin{lemma}  \label{l2.7}
Let $M$ and  $N$ be cellular pseudomanifolds, $\alpha \in M$ and
$\beta\in N$. \begin{enumerate}
     \item[$(a)$] $[({\bf 0}_{M}, \beta), ({\bf 1}_{M},
     \beta)]_{\otimes}$ is isomorphic to $M$ and $[(\alpha,
     {\bf 0}_{N}), (\alpha, {\bf 1}_{N})]_{\otimes}$ is isomorphic
     to $N$.
     \item[$(b)$] If $\beta$ is a facet of $N$ then
     $[({\bf 0}_{M}, \beta), ({\bf 1}_{M}, {\bf 1}_{N})]_{\ast}$
     is isomorphic to $M$ and if $\alpha$ is a facet of $M$ then
     $[(\alpha, {\bf 0}_{N}), ({\bf 1}_{M}, {\bf
     1}_{N})]_{\ast}$ is isomorphic to $N$.
\end{enumerate}
\end{lemma}

\noindent {\em Proof.} $\sigma \mapsto (\sigma, \beta)$ gives an
isomorphism between $M$ and $[({\bf 0}_{M}, \beta), ({\bf 1}_{M},
\beta)]_{\otimes}$. Similarly, $\gamma \mapsto (\alpha, \gamma)$
gives an isomorphism between $N$ and $[(\alpha, {\bf 0}_{N}),
(\alpha, {\bf 1}_{N})]_{\otimes}$. This proves $(a)$.

Let $\varphi \colon M \to [({\bf 0}_{M}, \beta), ({\bf 1}_{M},
{\bf 1}_{N})]_{\ast}$ be given by $\varphi(\sigma) = (\sigma,
\beta)$ for $\sigma < {\bf 1}_M$ and $\varphi({\bf 1}_M) = ({\bf
1}_M, {\bf 1}_N)$. Clearly, $\varphi$ is injective and $\tau <
\sigma$ implies that $\varphi(\tau) < \varphi(\sigma)$. If $(x,
y)\in [({\bf 0}_{M}, \beta), ({\bf 1}_{M}, {\bf 1}_{N})]_{\ast}$
then either $(x, y) = ({\bf 1}_M, {\bf 1}_N)$ or $y=\beta$. In the
first case, $\varphi({\bf 1}_M) = (x, y)$ and in the second case,
$\varphi(x) = (x, y)$. So, $\varphi$ is surjective and hence an
isomorphism. Similarly, $\psi \colon N \to [(\alpha, {\bf 0}_{N}),
({\bf 1}_{M}, {\bf 1}_{N})]_{\ast}$ given by $\varphi(\gamma) =
(\alpha, \gamma)$ for $\gamma < {\bf 1}_N$ and $\varphi({\bf 1}_N)
= ({\bf 1}_M, {\bf 1}_N)$ is an isomorphism. These prove $(b)$.
\hfill $\Box$

\begin{theo}  \label{t2.1} Let $M$ and $N$ be two
cellular pseudomanifolds.
\begin{enumerate}
    \item[$(i)$] For proper faces $\alpha$ and $\beta$ of
    $M$ and $N$ respectively, if $\partial \alpha$ and $\partial
    \beta$ are cellular spheres then the geometric carrier of the
    face $(\alpha, \beta) \in M \times N$ is pl-homeomorphic to
    $|\alpha| \times |\beta|$.
    In consequence, if the boundaries of all the
    facets in $M$ and $N$ are cellular spheres then $|M \times
    N|$ is pl-homeomorphic to $|M| \times |N|$.
    \item[$(ii)$] If $M$ and $N$ are cellular
    manifolds then $|M \times N|$ is pl-homeomorphic to $|M|
    \times |N|$ and hence $M \times N$ is a cellular manifold.
    \end{enumerate}
\end{theo}

\noindent{\em Proof.} Observe that if $f$ is a face of a cellular
pseudomanifold $K$ then $B(f)$ ($ = \{\sigma\cap [{\bf 0}, f] ~
\colon ~ \sigma\in B(M)\}$) is a simplicial complex and any
maximal simplex of $B(f)$ is of the form $\{f\}\cup \alpha$, where
$\alpha$ is a maximal simplex of $B(\partial f)$. This implies
that $|f|$ is homeomorphic to the cone over $|\partial f|$. (See
\cite{bd2} for the geometric carrier of a general simplicial
complex.)

If $f = (\alpha, \beta)$ then, by Lemma \ref{l2.6}, the set of
facets of $\partial f$ is $\{(g, \beta) ~ \colon ~ g \mbox{ is a
facet of } \partial \alpha\} \cup \{(\alpha, h) ~ \colon ~ h
\mbox{ is a facet of } \partial \beta\}$. So, $|\partial f|$ is
the boundary of the pl-ball $|\alpha| \times |\beta|$. Since a
pl-ball is pl-homeomorphic to the cone over its boundary (cf.
\cite{rs}), $|\alpha| \times |\beta|$ is pl-homeomorphic to the
cone over $|\partial f|$.  Since $|f|$ is the cone over $|\partial
f|$, $|\alpha| \times |\beta|$ is pl-homeomorphic to $|f|$. This
proves the first part of $(i)$. The last part of $(i)$ follows
from the first part.

Let $M$ and $N$ be cellular manifolds. If $\alpha$ is a facet of
$M$ then, by Lemma \ref{l2.3}, $\partial_M \alpha$ is a cellular
sphere. Similarly, $\partial_N \beta$ is a cellular sphere for
each facet $\beta$ of $N$. Therefore, by (i), $|M \times N|$ is
pl-homeomorphic to $|M| \times |N|$. Thus, $|M\times N|$ ($=
|B(M\times N)|$) is a pl-manifold and hence (cf. \cite{rs})
$B(M\times N)$ is a combinatorial manifold. This completes the
proof. \hfill $\Box$

\bigskip

If $M$ is a cellular pseudomanifold then, from the definition, $M
\ast S^{\,-1}_1$ is isomorphic to $M$. Here we prove\,:

\begin{theo}  \label{t2.2}
Let $M$ and $N$ be two cellular pseudomanifolds of dimensions
$\geq 0$.
\begin{enumerate}
    \item[$(i)$] The normal pseudomanifold $B(M\ast N)$ is a
    subdivision of $B(M)\ast B(N)$.
    \item[$(ii)$] If $M$ and $N$ are cellular spheres then $M
    \ast N$ is a cellular sphere.
    \item[$(iii)$] If $M$ and $N$ are polytopal spheres then $M
    \ast N$ is a polytopal sphere.
    \item[$(iv)$] If $M \ast N$ is a cellular manifolds then $M$
    and $N$ are cellular spheres and hence $M\ast N$ is a
    cellular sphere.
    \end{enumerate}
\end{theo}

\noindent {\em Proof.} Let $\dim(M) = m$ and $\dim(N) = n$.

If $\sigma$ is a proper face of a normal pseudomanifold $X$
and $v \not\in V(X)$ then $X^{\,\prime} := \{\alpha\in X ~ \colon
~ \sigma \not \subseteq \alpha\} \cup \{(\tau \setminus \{ x \})
\cup \{v\} ~ \colon ~ x \in \sigma \subseteq \tau\}$ is a
subdivision of $X$. We say that $X^{\,\prime}$ is obtained from
$X$ by {\em starring} a vertex in $\si$ (cf. \cite{bd2}).

Let $K_0 = B(M) \ast B(N)$. For $0 \leq i \leq m+ n$, let $E_i :=
\{(\alpha, \beta) ~ \colon ~ \alpha$ is a $j$-face of $M$ and
$\beta$ is a $(i-j)$-face of $N$, where $0\leq j \leq i\}
\subseteq$ the set of edges of $K_0$. Observe that any two edges
in $E_i$ are not in a 2-face (and hence not in any face) of $K_0$.
For $1\leq i\leq m + n + 1$, let $K_i$ be the subdivision of
$K_{i-1}$ obtained by starring a vertex on every edge in $E_{m + n
+1 -i}$. Since no two edges in $E_{m + n +1 -i}$ are in a face of
$K_0$ and hence in a face of $K_{i-1}$, it does not matter in
which ordering we do the starring on the edges in $E_{m + n +1
-i}$. Now, it is routine to check that $K_{m + n + 1}$ is
isomorphic to $B(M \ast N)$. This proves $(i)$.

Let $M$ and $N$ be cellular spheres. Since the join of two
combinatorial spheres is a combinatorial sphere (cf. \cite{bd2}),
$B(M)\ast B(N)$ is a combinatorial sphere. Then, by part $(i)$,
$B(M\ast N)$ is a combinatorial spheres. This proves $(ii)$.

Let $M$ and $N$ be the boundary complex of the polytopes $X
\subseteq \mathbb{R}^{m + 1}$ and $Y \subseteq \mathbb{R}^{n + 1}$ respectively.
Assume that $(0, \dots, 0)$ is in the interior of $X$ and $(0,
\dots, 0)$ is in the interior of $Y$. For $U \subseteq \mathbb{R}^{m + n
+ 2}$, let $\langle U \rangle$ denote the convex span of $U$. Let
$Z := \langle (X \times \{(0, \dots, 0)\}) \cup (\{(0, \dots, 0)\}
\times Y) \rangle$. Then $Z$ is a $(m + n + 2)$-polytope. Since
$(0, \dots, 0)$ is in the interior of $X$, $Y$ is not a closed
face of $Z$. Similarly, $X$ is not a closed face of $Z$. For a
hyperplane $P \subseteq \mathbb{R}^{m + n + 2}$ which touches $Z$, let
$P_1 := P \cap (\mathbb{R}^{m + 1} \times \{(0, \dots, 0)\})$ and $P_2 :=
P \cap (\{(0, \dots, 0)\} \times \mathbb{R}^{n + 1})$. Since $(0, \dots,
0)$ is in the interior of $Z$, $P_1 \neq \mathbb{R}^{m + 1} \times \{(0,
\dots, 0)\}$ and $P_2 \neq \{(0, \dots, 0)\} \times \mathbb{R}^{n + 1}$.
Therefore, $P \cap Z$ is either $A \times \{(0, \dots, 0)\}$,
$\{(0, \dots, 0)\} \times B$ or $\langle (A \times \{(0, \dots,
0)\}) \cup (\{(0, \dots, 0)\} \times B) \rangle$ for $A \neq X$ in
$M$, $B \neq Y$ in $N$. This implies that $\partial Z$ is
isomorphic to $M \ast N$. This proves $(iii)$.

Let $\beta$ be a facet of $N$. Then $({\bf 0}_M, \beta)$ is a
proper face of $M \ast N$ and hence, by Lemma \ref{l2.3}, $[({\bf
0}_M, \beta), ({\bf 1}_M, {\bf 1}_N)]$ is a cellular sphere. This
implies, by Lemma \ref{l2.7} $(b)$, that $M$ is a cellular sphere.
Similarly, $N$ is a cellular sphere. $(iv)$ now follows from part
$(ii)$. \hfill $\Box$

\begin{lemma}  \label{l2.8}
Let $M$ be a cellular pseudomanifold of dimension $d \geq 0$.
\begin{enumerate}
     \item[$(a)$] $S^{\,b}_{b + 2} \otimes S^{\,- 1}_1 = S^{\,b +
     1}_{b + 3}$ for $b \geq -1$. In consequence, the direct
     product of $k$ copies of $S^{\,- 1}_1$ is $= S^{\,k - 2}_k$ for
     $k\geq 2$ and $S^{\,b}_{b + 2} \otimes S^{\,c}_{c + 2} =
     S^{\,b + c + 2}_{b + c + 4}$ for $b, c \geq -1$.
     \item[$(b)$] $M \otimes S^{\,-1}_1$ is a cellular
     pseudomanifold of dimension $d + 1$ and $e(M \otimes
     S^{\,-1}_1) = e(M)$.
     \item[$(c)$] $B(M \otimes S^{\,-1}_1)$ is isomorphic to a
     subdivision of $B(M) \ast S^{\,0}_2$. In consequence, the
     geometric carrier of $M \otimes S^{\,-1}_1$ is homeomorphic
     to the suspension of the geometric carrier of $M$. Hence the
     geometric carrier of $M \otimes S^{\,c}_{c + 2}$ is the $(c
     + 2)$-fold suspension of the geometric carrier of $M$.
     \item[$(d)$] $M$ is a cellular sphere if and only if $M \otimes
     S^{\,-1}_1$ is a cellular sphere. In consequence, $M$ is a cellular
     sphere if and only if $M \otimes S^{\,c}_{c + 2}$ is a cellular
     sphere for $c \geq -1$.
     \item[$(e)$] If $M$ is a polytopal sphere then $M\otimes S^{\,-
     1}_1$ is a polytopal sphere. In consequence, if $M$ is a polytopal
     sphere then $M \otimes S^{\,c}_{c + 2}$ is a polytopal sphere for
     $c \geq -1$.
     \item[$(f)$] If $M \otimes S^{\,-1}_{1} = N \otimes S^{\,-1}_{1}$
     then $M = N$. In consequence, if $M \otimes S^{\,c}_{c + 2} = N
     \otimes S^{\,c}_{c + 2}$ then $M = N$.
     \end{enumerate}
\end{lemma}

\noindent {\em Proof.} $(a)$, $(b)$ and $(f)$ follow from the
definition of direct product.

To prove $(c)$, let $K_0 = B(M) \ast S^{\,0}_2(\{v, w\})$. For $0
\leq i \leq d$, let $E_i := \{\{w, \sigma \} ~ \colon ~ \sigma$
is an $i$-face of $M\} \subseteq$ the set of edges of $K_0$.
Observe that any two edges in $E_i$ are not in a 2-face (and
hence not in any face) of $K_0$. For $1 \leq i\leq d + 1$, let
$K_i$ be the subdivision of $K_{i-1}$ obtained by starring a
vertex on every edge in $E_{d + 1 -i}$. Since no two edges in
$E_{d + 1 -i}$ are in a face of $K_0$ and hence in a face of
$K_{i - 1}$, it does not matter in which ordering we do the
starring on the edges in $E_{d + 1 -i}$. Let the new vertex in
the edge $\{w, \sigma\}$ be $v_{\sigma}$. So, $K_{d+1}$ is a
subdivision of $B(M) \ast S^{\,0}_2(\{v, w\})$.

Consider the bijection $\varphi \colon V(B(M \otimes
S^{-1}_1(\{u\}))) \to V(K_{d+1})$ given by $\varphi({\bf 1}_M) =
v$, $\varphi(u) = w$, $\varphi(\tau) = \tau$ for $\tau \in
V(B(M))$ and $\varphi(\sigma \cup \{w\}) = v_{\sigma}$ for
$\sigma \in M$. Now, it is routine to check that $\varphi$ is an
isomorphism. This proves $(c)$.

If $M$ is a cellular sphere then $B(M)$ is a combinatorial sphere
and hence $B(M) \ast S^{\,0}_2$ is a combinatorial sphere. Then,
by $(c)$, $B(M \otimes S^{\,-1}_1)$ is a combinatorial sphere and
hence $M \otimes S^{\,-1}_1$ is a cellular sphere. Conversely, if
$M \otimes S^{\,-1}_1$ a cellular sphere then $B(M \otimes
S^{\,-1}_1)$ is a combinatorial sphere. Therefore, by $(c)$, $N
:= B(M) \ast S^{\,0}_2(\{u, v\})$ is a combinatorial sphere and
hence a combinatorial manifold. Therefore, since ${\rm lk}_N(u) =
B(M)$, $B(M)$ is a combinatorial sphere and hence $M$ is a
cellular sphere. This proves $(d)$.

Let $X \subseteq \mathbb{R}^{d + 1}$ be an $(d + 1)$-polytope such that
$\partial X = M$. (We may assume that $(0, \dots, 0)$ is an
interior point of $X$.) Let $Y$ be the convex span of $(X \times
\{0\}) \cup \{(0, \dots, 0, 1)\}$ in $\mathbb{R}^{d + 2}$. Then $Y$ is a
$(d + 2)$-polytope and the boundary complex of $Y$ is isomorphic
to $M \otimes S^{\,-1}_1(\{v\})$ (under the map given by $(u, 0)
\mapsto (u, \emptyset) \equiv u$ for $u \in V(M)$ and $(0, \dots,
0, 1) \mapsto ({\bf 0}_M, v) \equiv v$). This proves $(e)$.
\hfill $\Box$

\begin{theo}  \label{t2.3}
Let $M$ and $N$ be two cellular pseudomanifolds.
\begin{enumerate}
    \item[$(i)$] If $M$ and $N$ are cellular spheres then $M
    \otimes N$ is a cellular sphere.
    \item[$(ii)$] If $M$ and $N$ are polytopal spheres then $M
    \otimes N$ is a polytopal sphere.
    \item[$(iii)$] If $M \otimes N$ is a cellular manifolds then
    $M$ and $N$ are cellular spheres and hence $M \otimes N$ is
    a cellular sphere.
    \end{enumerate}
\end{theo}

\noindent {\em Proof.} Observe that $M^{\,\prime} := \{(\alpha,
\emptyset) \, : \, \alpha \in M\} \subseteq M \otimes S^{\,-1}_1$
is isomorphic to $M$ by the map $(\alpha, \emptyset) \mapsto
\alpha$ and $N^{\,\prime} := \{(\beta, \emptyset) \, : \, \beta
\in N\} \subseteq N \otimes S^{\,-1}_1$ is isomorphic to $N$ by
the map $(\beta, \emptyset) \mapsto \beta$. (Here $\emptyset$ is
the ${\bf 0}$ of $S^{\,-1}_1$.)

If $M$, $N$ are cellular spheres then, by Lemma \ref{l2.8} $(d)$,
$M \otimes S^{\,-1}_1$ and $N \otimes S^{\,-1}_1$ are cellular
spheres. Therefore, by Theorem \ref{t2.2} $(ii)$, $P := (M \otimes
S^{\,-1}_1) \ast (N \otimes S^{\,-1}_1)$ is a cellular sphere and
hence a cellular manifold. Since $({\bf 1}_M, \emptyset)$ is a
proper face of $M \otimes S^{\,-1}_1$ and $({\bf 1}_N, \emptyset)$
is a proper face of $N \otimes S^{\,-1}_1$, $\alpha := (({\bf
1}_M, \emptyset), ({\bf 1}_N, \emptyset))$ is a proper face (in
fact, a facet) of $P$. Therefore, by Lemma \ref{l2.3}, $\partial_P
\alpha$ is a cellular sphere. Now, if $\beta \in \partial_P
\alpha$ then $\beta = ((f, \emptyset), (g, \emptyset))$ for some
$f \in [{\bf 0}_M, {\bf 1}_M]$ and $g \in [{\bf 0}_N, {\bf 1}_N]$.
So, $\partial_P \alpha$ is isomorphic to $M^{\,\prime} \otimes
N^{\,\prime}$. This implies that $M \otimes N$ is isomorphic to
$\partial_P \alpha$. This proves $(i)$.

Let $\dim(M) = m - 1$ and $\dim(N) = n - 1$. Let $M$ and $N$ be
the boundary complex of the polytopes $X \subseteq \mathbb{R}^{m}$ and $Y
\subseteq \mathbb{R}^{n}$ respectively. Then $X \ast Y$ is an $(m + n +
1)$-polytope and, by Lemma \ref{l2.4}, $M \ast N$ is isomorphic to
$\partial (X\ast Y)$. This proves $(ii)$.

Let $\beta$ be a face of $N$. Then $[({\bf 0}_M, \beta), ({\bf
1}_M, \beta)]\neq M \otimes N$ and hence, by Lemma \ref{l2.3},
$[({\bf 0}_M, \beta), ({\bf 1}_M, \beta)]$ is a cellular sphere.
This implies, by Lemma \ref{l2.7} $(a)$, that $M$ is a cellular
sphere. Similarly, $N$ is a cellular sphere. $(iii)$ now follows
from part $(i)$. \hfill $\Box$

\section{\sc Cellular pseudomanifolds of excess one}

In \cite[Section 6.1]{g}, Gr\"{u}nbaum has classified all the
combinatorial types of $(d+1)$-polytopes with $d+3$ vertices.
There are exactly $\lfloor(\frac{d+1}{2})^2 \rfloor$ such
combinatorial types. This amounts to a classification of the
polytopal spheres of dimension $d$ and excess 1. In this section,
we present a classification of all the $d$-dimensional cellular
pseudomanifolds with excess 1. It turns out that all of them are
polytopal spheres. Thus\,:

\begin{theo}  \label{t3.1}
The only $d$-dimensional cellular pseudomanifolds of excess $1$
are $S^{\,b_1}_{b_1 +2} \ast S^{\,b_2}_{b_2 +2}$ with $0\leq b_1
\leq b_2$, $b_1 + b_2 = d-1$ and $(S^{\,d_1}_{d_1 +2} \ast
S^{\,d_2}_{d_2 +2}) \otimes S^{\,d_3}_{d_3 +2}$ with $0\leq d_1
\leq d_2$, $-1 \leq d_3$, $d_1 + d_2 + d_3 = d-3$. These are
mutually non-isomorphic polytopal spheres. There are exactly
$\lfloor (\frac{d + 1}{2})^2\rfloor$ of them.
\end{theo}

\noindent {\em Proof.} Let $M$ be a $d$-dimensional cellular
pseudomanifold of excess 1, i.e., with $d+3$ vertices. As before,
we identify a face with its shadow. If each facet contains $d+1$
vertices then $M$ is simplicial. In that case, by the
classification of pseudomanifolds of excess 1 (in
\cite{bd2}), $M= S^{\,b_1}_{b_1+2} \ast S^{\,b_2}_{b_2 + 2}$ for
some $b_1$, $b_2$ with $0\leq b_1 \leq b_2$, $b_1 + b_2 = d - 1$.

So, assume that there exists a facet $\sigma$ containing $d+2$
vertices. Then $(\partial \sigma) \otimes S^{\,-1}_{1}(\{v\})
\subseteq M$, where $v$ is the unique vertex outside $\sigma$.
(Here we identify $(A, B)$ with $A\cup B$.) By connectedness of
$\Lambda(M)$, $(\partial \sigma) \otimes S^{\,-1}_{1}(\{v\}) = M$.

Let $d_3\geq -1$ be the largest integer such that $M = N \otimes
S^{\,d_3}_{d_3 + 2}$ for some cellular pseudomanifold $N$. Then,
by Lemma \ref{l2.5}, $e(N) = 1$ and $N$ is simplicial (else $N =
N_0 \otimes S^{\,-1}_1$ for some cellular pseudomanifold $N_0$ and
hence, by Lemma \ref{l2.8} $(a)$, $M = N_0 \otimes S^{\,d_3 +
1}_{d_3 +3}$, contradicting the choice of $d_3$). Therefore, $N =
S^{\,d_1}_{d_1+2} \ast S^{\,d_2}_{d_2 + 2}$ for some $d_1$, $d_2$
with $0\leq d_1 \leq d_2$, $d_1 + d_2 = d - d_3- 3$. This proves
the first part.

Observe that $S^{\,b_1}_{b_1 +2} \ast S^{\,b_2}_{b_2 +2}$ is
simplicial but $(S^{\,d_1}_{d_1 +2} \ast S^{\,d_2}_{d_2 +2})
\otimes S^{\,d_3}_{d_3 +2}$ is not simplicial. Therefore, by Lemma
\ref{l2.8} $(f)$, all the cellular pseudomanifolds in the
statement are mutually non-isomorphic. Since a normal 
pseudomanifold of type $S^{\,b_1}_{b_1+2} \ast S^{\,b_2}_{b_2 +
2}$ is a polytopal sphere, by Lemma \ref{l2.8} $(e)$, all the
cellular pseudomanifolds in the statement are polytopal spheres.

Since there are $\lfloor\frac{d+1}{2}\rfloor$ choice of $b_1$ and
$b_2$ is determined by $b_1$, the total number of cellular
pseudomanifolds of type $S^{\,b_1}_{b_1 +2} \ast S^{\,b_2}_{b_2
+2}$ is $\lfloor\frac{d+1}{2}\rfloor$.

Thus, for each $d_3$ ($- 1 \leq d_3 \leq d - 3$), there are
$\lfloor \frac{d - d_3 - 1}{2} \rfloor$ cellular pseudomanifolds
of type $(S^{\,d_1}_{d_1 + 2} \ast S^{\,d_2}_{d_2 + 2}) \otimes
S^{\,d_3}_{d_3 + 2}$. Hence the total number of cellular
pseudomanifolds of second type is $\sum_{d_3 = -1}^{d -3} \lfloor
\frac{d - d_3 - 1}{2} \rfloor = \lfloor\frac{d^2}{4}\rfloor$.
Adding, we get the total number of cellular pseudomanifolds of
excess 1 is $ \lfloor \frac{d^2}{4} \rfloor + \lfloor \frac{d +
1}{2}\rfloor = \lfloor (\frac{d + 1}{2})^2 \rfloor$. \hfill $\Box$

\section{\sc Equivalence relation $\sim$}

\begin{defn}  {\rm Let $M$ be a cellular pseudomanifold. Define
the binary relation `$\sim$' on $V(M)$ as follows. For $x, y\in
V(M)$, $x\sim y$ if and only if either $x = y$ or the
transposition $(x, y)$ interchanging $x$ and $y$ (and fixing all
other vertices) is an automorphism of $M$. It is easy to verify
that $\sim$ is an equivalence relation. (Whenever $x$, $y$, $z$
are distinct vertices with $x\sim y$ and $y\sim z$, we have $(x,
z) = (x, y)(y, z)(x, y)$ is an automorphism and hence $x\sim z$.)}
\end{defn}

\begin{lemma}  \label{l4.1} Let $M$ be a cellular pseudomanifold
and let $A$ be a $\sim$-equivalence class.
\begin{enumerate}
    \item[$(a)$]  If $\beta$ is a face such that $A\not\subseteq
    {\rm sh}(\beta)$ then ${\rm sh}(\beta) \setminus A$ is the
    shadow of a face.
    \item[$(b)$] If $\sigma$ is the smallest face such that ${\rm
    sh}(\sigma)$ contains $A$ then ${\rm sh}(\sigma)$ is a union
    of $\sim$-equivalence classes.
    \end{enumerate}
\end{lemma}

\noindent {\em Proof.} In this proof, we identify a face with its
shadow. If $\beta\cap A = \emptyset$ then there is nothing to
prove. Fix a vertex $v\in A\setminus \beta$. Then $\beta_x := (v,
x)(\beta)$ is a face for each $x\in A\cap \beta$. Since $x\not\in
\beta_x$, $\beta\cap(\cap_{x\in A\cap \beta}\beta_x) =
\beta\setminus A$. Hence $\beta \setminus A$ is a face. This
proves $(a)$.

If there exists a $\sim$-equivalence class $B$ such that
$B\cap\sigma\neq\emptyset$ and $B\not\subseteq \sigma$ then, by
part $(a)$, $\sigma\setminus B$ is a face containing $A$. This
contradicts the minimality of $\sigma$. This proves $(b)$. \hfill
$\Box$

\begin{lemma}  \label{l4.2} Let $M$ be a cellular
pseudomanifold.  \begin{enumerate}
    \item[$(a)$]  For $x, y \in V(M)$, $x \sim y$ implies that
    every facet of $M$ is incident with either $x$ or $y$.
    \item[$(b)$] If $A$ is a $\sim$-equivalence class and
    $\alpha$ is a face such that ${\rm sh}(\alpha)\cap A =
    \emptyset$ then for any proper subset $B$ of $A$, ${\rm
    sh}(\alpha) \cup B$ is the shadow of a face.
    \item[$(c)$] More generally, for $U\subseteq V(M)$ if
    $\widehat{U}$ denotes the union of all the $\sim$-equivalence
    classes contained in $U$ then $U$ is the shadow of a face
    if and only if $\widehat{U}$ is the shadow of a face.
    \item[$(d)$] Let $\sigma$ and $\tau$ be two facets adjacent
    in $\Lambda(M)$. Let $y$ be a vertex incident with $\sigma$
    but not with $\tau$. If $x\sim y$ then $(x, y)(\tau)$ is
    either equal to $\sigma$ or adjacent with $\sigma$ in
    $\Lambda(M)$.
    \end{enumerate}
\end{lemma}

\noindent {\em Proof.} Here, we identify a face with its shadow.
The proof of $(a)$ is by induction on $\dim(M)$. The result is
obvious for dimension $\leq 1$. So, let $\dim(M) > 1$ and suppose
that we have the result for all cellular pseudomanifolds of
dimension $< \dim(M)$.

Let ${\cal C}$ be the set of all facets which are incident with
either $x$ or $y$. Note that ${\cal C}$ is a non-empty set of
vertices of the graph $\Lambda(M)$.  Let $f \in {\cal C}$ and let
$g$ be a neighbour of $f$ in $\Lambda (M)$. First assume $x$, $y
\in f$. Then $(x, y)$ is an automorphism of the cellular
pseudomanifold $\partial f$. If $x \not\in g$, $y \not\in g$ then
$f\wedge g$ is a facet of $\partial f$ not containing $x$ or $y$
and $\dim(\partial f)= \dim(M)-1$. Contradiction to the induction
hypothesis. So, at least one of $x$, $y$ is in $g$ and hence $g\in
{\cal C}$. Next suppose $f$ contains exactly one of $x$, $y$. Say
$x \in f$, $y \not\in f$. If $x \in g$ then $g \in {\cal C}$.
Otherwise, the image of $f$ under the automorphism $(x, y)$
contains $f \wedge g$ and is $\neq f$ and hence the image is $g$.
Thus $y\in g$, and hence $g \in {\cal C}$. Thus we see that in all
cases, the neighbours in $\Lambda(M)$ of each element of ${\cal
C}$ is again in ${\cal C}$. By connectedness of $\Lambda(M)$,
$\cal C$ consists of all the facets of $M$. This proves $(a)$.

By Lemma \ref{l2.1}, for each $x\in V(M)\setminus(\alpha\cup B)$,
there exists a facet $\sigma_x$ such that $\alpha\subseteq
\sigma_x$, $x\not\in\sigma_x$. By part $(a)$, $\sigma_x$ misses at
most one vertex from $A$. Applying a suitable automorphism to
$\sigma_x$ (if necessary), we may assume that $\sigma_x$ contains
$B$. Then $\alpha\cup B = \cap_x \sigma_x$ is a face. This proves
$(b)$.

To prove $(c)$, let $B_1, \dots, B_k$ be all the
$\sim$-equivalence classes such that $U\cap B_i\neq\emptyset$ and
$B_i\not\subseteq U$ for $1\leq i\leq k$. Now, if $U$ is a face
then, by part $(a)$ of Lemma \ref{l4.1}, inductively we get that
each of $U \setminus B_1$, $(U\setminus B_{1})\setminus B_2,
\dots, (U \setminus (B_{1} \cup \cdots \cup B_{k-1}))\setminus B_k
= \widehat{U}$ is a face. Conversely, if $\widehat{U}$ is a face
then put $C_i := B_i\cap U$ for $1\leq i\leq k$. Using part $(b)$,
we get inductively $\widehat{U} \cup C_1$, $(\widehat{U}\cup C_1)
\cup C_2, \dots, (\widehat{U} \cup C_1\cup \cdots \cup C_{k-1})
\cup C_k= U$ are faces. This proves $(c)$.

To prove $(d)$, let $C$ be the $\sim$-equivalence class of $x$ and
$y$. Since $\tau$ is a facet and $y\not\in \tau$, $x\in \tau$ and
$C\not\subseteq \tau$. This implies, by part $(a)$ of Lemma
\ref{l4.1}, $\tau \setminus C$ is a face. Then, by part $(b)$,
$\tau \setminus \{x\} = (\tau \setminus C)\cup ((\tau\cap C)
\setminus \{x\})$ is a face.

If $x\not\in\sigma$ then, by the same argument as above, $y\in
\sigma$ and $\sigma\setminus\{y\}$ is a face. Now, $x, y \not \in
\sigma \cap \tau$. Therefore, $\sigma \cap \tau \subseteq \tau
\setminus \{x\} \subset \tau$ and $\sigma \cap \tau \subseteq
\sigma \setminus \{y\} \subset \sigma$. Since $\dim(\sigma \cap
\tau) = \dim(\tau) - 1 = \dim(\sigma) - 1$, $\sigma \cap \tau =
\tau \setminus \{x\} = \sigma \setminus \{y\}$. This implies that
$\sigma = (\tau \setminus \{x\}) \cup \{y\} = (x, y)(\tau)$.

If $x\in \sigma$ then $\sigma \cap ((x, y)(\tau)) = ((x,
y)(\sigma)) \cap ((x, y)(\tau)) = (x, y)(\sigma\cap\tau)$. Since
$(x, y)$ is an automorphism, $(x, y)(\sigma\cap\tau)$ is a face of
dimension $\dim(\sigma)-1$. These imply that $\sigma$ and $(x,
y)(\tau)$ are adjacent in $\Lambda(M)$. This proves $(d)$. \hfill
$\Box$

\begin{lemma}  \label{l4.3} Let $M$ be a cellular
pseudomanifold. \begin{enumerate}
    \item[$(a)$] Let $f\subseteq V(M)$ be a set of vertices of $M$
    such that $f$ misses exactly one vertex from each
    $\sim$-equivalence class. Then $f$ is $($the shadow of$\,)$ a
    face of dimension $\#(f)-1$.
    \item[$(b)$] Let $f$, $f^{\,\prime}$ be two faces of $M$, as
    in part $(a)$. Then there is an automorphism $\varphi$ of $M$
    such that $f^{\,\prime} = \varphi(f)$.
\end{enumerate}
\end{lemma}

\noindent {\em Proof.} As usual, we identify a face with its
shadow.

By part $(c)$ of Lemma \ref{l4.2}, any subset of $f$ is a face.
This implies that $f$ is a face of rank $\#(f)$. This proves
$(a)$.

Let $A_1, \dots, A_m$ be the $\sim$-equivalence classes. Put
$\{x_i\} = A_i\setminus f$ and  $\{x_i^{\,\prime}\} = A_i\setminus
f^{\,\prime}$ for $1\leq i\leq m$. Take $\varphi$ to be the
product of those transpositions $(x_i, x_i^{\,\prime})$ where
$x^{\prime}_i \neq x_i$. Then $f^{\,\prime} = \varphi(f)$. This
proves $(b)$. \hfill $\Box$

\begin{eg}  \label{v6d2} {\rm Some 2-dimensional
cellular pseudomanifolds on 5 and 6 vertices. }
\end{eg}

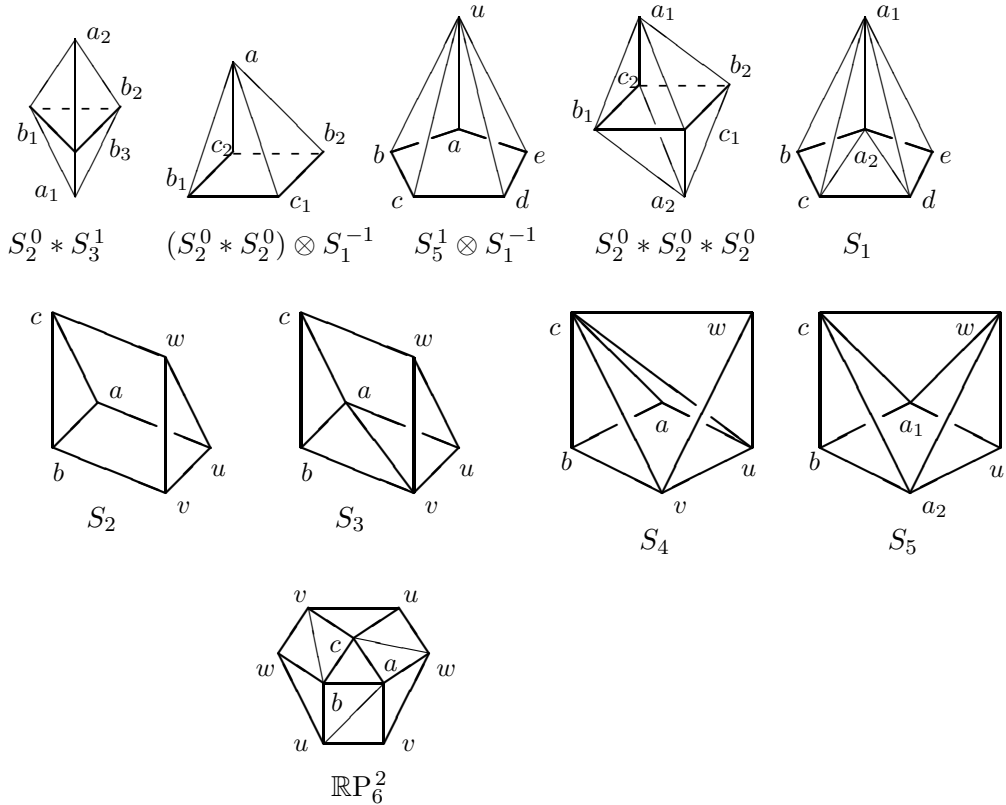
\setlength{\unitlength}{3mm}

\begin{picture}(48.5,12)(1,0)


\thicklines

\put(4,5){\line(-1,1){2}} \put(4,5){\line(1,1){2}}

\thinlines

\put(4,3){\line(-1,2){2}} \put(4,3){\line(1,2){2}}
\put(4,3){\line(0,1){7}} \put(4,10){\line(-2,-3){2}}
\put(4,10){\line(2,-3){2}}

\put(2.55,6.65){- - - -}

\put(2.2,3){\mbox{\small $a_{1}$}} \put(4.5,10){\mbox{\small
$a_{2}$}} \put(1.4,5.5){\small \mbox{$b_{1}$}}
\put(6,7.5){\mbox{\small $b_{2}$}} \put(5.5,4.8){\mbox{\small
$b_{3}$}}

\put(1,0.5){\mbox{$S^{\,0}_{2}\ast S^{\,1}_3$}}

\thicklines

\put(9,3){\line(1,0){4}} \put(9,3){\line(1,1){2}}
\put(13,3){\line(1,1){2}}

\thinlines

\put(9,3){\line(1,3){2}} \put(11,5){\line(0,1){4}}
\put(13,3){\line(-1,3){2}} \put(15,5){\line(-1,1){4}}

\put(11,4.65){- - - - -}

\put(11.5,9){\mbox{\small $a$}} \put(8,3.5){\mbox{\small $b_{1}$}}
\put(13.5,2.5){\small \mbox{$c_{1}$}} \put(15,5.5){\mbox{\small
$b_{2}$}} \put(10,5){\mbox{\small $c_{2}$}}

\put(8,0.5){\mbox{$(S^{\,0}_{2}\ast S^{\,0}_2)\otimes S^{\,-1}_{1}
$}}


\thicklines

\put(19,3){\line(-1,2){1}} \put(19,3){\line(1,0){4}}
\put(23,3){\line(1,2){1}} \put(18,5){\line(3,1){1.2}}
\put(21,6){\line(3,-1){1.18}} \put(24,5){\line(-3,1){1.2}}
\put(21,6){\line(-3,-1){1.18}}

\thinlines

\put(18,5){\line(1,2){3}} \put(19,3){\line(1,4){2}}
\put(23,3){\line(-1,4){2}} \put(24,5){\line(-1,2){3}}
\put(21,6){\line(0,1){5}}

\put(21.5,11){\mbox{\small $u$}} \put(20.5,5){\mbox{\small $a$}}
\put(17.2,4.5){\mbox{\small $b$}} \put(24.3,4.5){\mbox{\small
$e$}} \put(18,2.5){\mbox{\small $c$}} \put(23.5,2.5){\mbox{\small
$d$}}

\put(19,0.5){\mbox{$S^{\,1}_{5}\otimes S^{\,-1}_{1}$}}


\thicklines

\put(27,6){\line(1,0){4}} \put(27,6){\line(1,1){2}}
\put(31,6){\line(1,1){2}}

\thinlines

\put(27,6){\line(2,5){2}} \put(29,8){\line(0,1){3}}
\put(31,6){\line(-2,5){2}} \put(33,8){\line(-4,3){4}}
\put(27,6){\line(4,-3){4}} \put(31,3){\line(2,5){2}}
\put(31,3){\line(0,1){3}} \put(31,3){\line(-2,5){1.05}}
\put(29,8){\line(2,-5){0.7}}

\put(29,7.65){- - - - -}

\put(29.5,11){\mbox{\small $a_1$}} \put(29.5,2.5){\small
\mbox{$a_{2}$}} \put(26,6.5){\mbox{\small $b_{1}$}}
\put(32.5,5.5){\small \mbox{$c_{1}$}} \put(33,8.5){\mbox{\small
$b_{2}$}} \put(28,8){\mbox{\small $c_{2}$}}

\put(27,0.5){\mbox{$S^{\,0}_{2}\ast S^{\,0}_{2}\ast S^{\,0}_{2}$}}


\thicklines

\put(37,3){\line(-1,2){1}} \put(37,3){\line(1,0){4}}
\put(41,3){\line(1,2){1}} \put(36,5){\line(3,1){1.2}}
\put(39,6){\line(3,-1){1.18}} \put(42,5){\line(-3,1){1.2}}
\put(39,6){\line(-3,-1){1.18}}

\thinlines

\put(36,5){\line(1,2){3}} \put(37,3){\line(1,4){2}}
\put(41,3){\line(-1,4){2}} \put(42,5){\line(-1,2){3}}
\put(39,6){\line(0,1){5}} \put(39,6){\line(-2,-3){2}}
\put(39,6){\line(2,-3){2}}

\put(39.5,11){\mbox{\small $a_1$}} \put(38.5,4.5){\mbox{\small
$a_2$}} \put(35.2,4.5){\mbox{\small $b$}}
\put(42.3,4.5){\mbox{\small $e$}} \put(36,2.5){\mbox{\small $c$}}
\put(41.5,2.5){\mbox{\small $d$}}

\put(38,0.5){\mbox{$S_1$}}

\end{picture}

\setlength{\unitlength}{3mm}

\begin{picture}(45,13)(0,0)


\thicklines

\put(2,5){\line(5,-2){5}} \put(2,5){\line(1,1){2}}
\put(2,5){\line(0,1){6}} \put(7,3){\line(1,1){2}}
\put(7,3){\line(0,1){6}} \put(9,5){\line(-1,2){2}}
\put(4,7){\line(-1,2){2}} \put(2,11){\line(5,-2){5}}
\put(4,7){\line(5,-2){2.6}} \put(9,5){\line(-5,2){1.6}}


\put(2,3.5){\mbox{\small $b$}} \put(1,10.5){\mbox{\small $c$}}
\put(4.5,7.2){\small \mbox{$a$}} \put(9,3.75){\mbox{\small $u$}}
\put(7.5,2.1){\mbox{\small $v$}} \put(7,9.5){\mbox{\small $w$}}

\put(3.5,1.5){\mbox{$S_2$}}


\thicklines

\put(13,5){\line(5,-2){5}} \put(13,5){\line(1,1){2}}
\put(13,5){\line(0,1){6}} \put(18,3){\line(1,1){2}}
\put(18,3){\line(0,1){6}} \put(20,5){\line(-1,2){2}}
\put(15,7){\line(-1,2){2}} \put(13,11){\line(5,-2){5}}
\put(15,7){\line(5,-2){2.6}} \put(20,5){\line(-5,2){1.6}}
\put(18,3){\line(-3,4){3}}


\put(13,3.5){\mbox{\small $b$}} \put(12,10.5){\mbox{\small $c$}}
\put(15.5,7.2){\small \mbox{$a$}} \put(20,3.75){\mbox{\small $u$}}
\put(18.5,2.1){\mbox{\small $v$}} \put(18,9.5){\mbox{\small $w$}}

\put(14.5,1.5){\mbox{$S_3$}}


\thicklines

\put(25,5){\line(2,-1){4}} \put(25,5){\line(0,1){6}}
\put(33,5){\line(-2,-1){4}} \put(33,5){\line(0,1){6}}
\put(29,3){\line(1,2){4}} \put(29,7){\line(-1,1){4}}
\put(25,11){\line(1,0){8}} \put(33,5){\line(-2,1){2}}
\put(29,7){\line(2,-1){1.2}}

\put(25,5){\line(2,1){2}} \put(29,7){\line(-2,-1){1.2}}
\put(29,3){\line(-1,2){4}} \put(33,5){\line(-4,3){1.8}}
\put(25,11){\line(4,-3){5.5}}

\put(24.52,4){\mbox{\small $b$}} \put(24,10){\mbox{\small $c$}}
\put(28.7,5.7){\small \mbox{$a$}} \put(32.5,3.75){\mbox{\small
$u$}} \put(29.5,2.1){\mbox{\small $v$}} \put(31,10){\mbox{\small
$w$}}

\put(28,0.5){\mbox{$S_4$}}


\thicklines

\put(36,5){\line(2,-1){4}} \put(36,5){\line(0,1){6}}
\put(44,5){\line(-2,-1){4}} \put(44,5){\line(0,1){6}}
\put(40,3){\line(1,2){4}} \put(40,7){\line(-1,1){4}}
\put(36,11){\line(1,0){8}} \put(44,5){\line(-2,1){2}}
\put(40,7){\line(2,-1){1.2}}

\put(36,5){\line(2,1){2}} \put(40,7){\line(-2,-1){1.2}}
\put(40,3){\line(-1,2){4}} \put(40,7){\line(1,1){4}}

\put(35.5,4){\mbox{\small $b$}} \put(35,10){\mbox{\small $c$}}
\put(39.5,5.7){\small \mbox{$a_{1}$}} \put(43.5,3.7){\mbox{\small
$u$}} \put(40.5,2.1){\mbox{\small $a_{2}$}}
\put(42,10){\mbox{\small $w$}}

\put(39,0.5){\mbox{$S_5$}}

\end{picture}

\setlength{\unitlength}{2mm}

\begin{picture}(46,16)(0,-1)


\thicklines

\put(21,3){\line(1,0){4}} \put(21,3){\line(0,1){4}}
\put(21,3){\line(-1,2){3}} \put(21,7){\line(2,3){2}}
\put(21,7){\line(-3,2){3}} \put(21,7){\line(1,0){4}}
\put(25,7){\line(-2,3){2}} \put(25,7){\line(3,2){3}}
\put(25,3){\line(1,2){3}} \put(25,3){\line(0,1){4}}
\put(23,10){\line(3,2){3}} \put(23,10){\line(-3,2){3}}
\put(18,9){\line(2,3){2}} \put(28,9){\line(-2,3){2}}
\put(20,12){\line(1,0){6}}

\thinlines

\put(21,3){\line(1,1){4}} \put(21,7){\line(-1,5){1}}
\put(28,9){\line(-5,1){5}}

\put(19,2.5){\mbox{\small $u$}} \put(26.2,2.5){\mbox{\small $v$}}
\put(28.5,7.5){\mbox{\small $w$}} \put(16.5,7.5){\mbox{\small
$w$}} \put(19,12.4){\mbox{\small $v$}}
\put(26.3,12.4){\mbox{\small $u$}} \put(25,7.8){\small \mbox{$a$}}
\put(21.4,9.1){\mbox{\small $c$}} \put(21.5,5.2){\mbox{\small
$b$}}

\put(21.5,-0.5){\mbox{$\mathbb{R}\mbox{P}^{\,2}_6$}}

\end{picture}

\medskip

The $\sim$-equivalence classes of $S^{\,1}_3 \ast S^{\,0}_2$ are
$\{a_1, a_2\}$ and $\{b_1, b_2, b_3\}$. The $\sim$-equivalence
classes of $(S^{\,0}_2\ast S^{\,0}_2)\otimes S^{\,-1}_1$ are
$\{a\}$, $\{b_1, b_2\}$ and $\{c_1, c_2\}$. The $\sim$-equivalence
classes of $S^{\,0}_2 \ast S^{\,0}_2 \ast S^{\,0}_2$ are $\{a_1,
a_2\}$, $\{b_1, b_2\}$ and $\{c_1, c_2\}$. The only non-trivial
$\sim$-equivalence class in $S_1$ is $\{a_1, a_2\}$. The
non-trivial $\sim$-equivalence class of $S_5$ is $\{a_1, a_2\}$.

\section{\sc Reducible cellular pseudomanifolds}

\begin{defn}  \label{reducible}
{\rm A cellular pseudomanifold $M$ is said to be} reducible {\rm
if there is a $\sim$-equivalence class of $M$ which is not
contained in the shadow of any facet of $M$. We shall say that $M$
is} irreducible {\em if $M$ is not reducible.}
\end{defn}

$S^{\,1}_3 \ast S^{\,0}_2$ and $S^{\,0}_2 \ast S^{\,0}_2 \ast
S^{\,0}_2$ are reducible but all others in Example \ref{v6d2} are
irreducible. The following characterization of reducibility in
terms of join decomposition justifies its name.

\begin{theo}  \label{t5.1} A cellular pseudomanifold
$M$ is reducible if and only if $M = N \ast S^{\,c}_{c + 2}$ for
some cellular pseudomanifold $N$ and integer $c \geq 0$. $($Notice
that, in particular, $S^{\,c}_{c + 2} = S^{\,- 1}_1 \ast
S^{\,c}_{c + 2}$ is reducible.$)$
\end{theo}

This theorem will be more or less immediate from the following\,:

\begin{lemma}  \label{l5.1} Let $A$ be a
$\sim$-equivalence class of a cellular pseudomanifold $M$ such
that $A$ is not contained in the shadow of any facet of $M$. Let
$N$ denote the subset of $M$ consisting of all faces whose shadows
are disjoint from $A$. Then $($after throwing in a trivial face
${\bf 1})$ $N$ is a cellular pseudomanifold.
\end{lemma}

\noindent {\em Proof.} Fix a vertex $a \in A$, and let $F =
A\setminus \{a\}$. By Lemma \ref{l4.2} $(b)$, $F\in M$ (with the
usual identification of face with its shadow). We claim that $N$
is isomorphic to ${\rm lk}_M(F)$, and hence a cellular
pseudomanifold. Indeed, by Lemmas \ref{l4.1} $(a)$ and \ref{l4.2}
$(b)$, the map $\alpha \mapsto \alpha\cup F$ for $\alpha < {\bf
1}_N$ and ${\bf 1}_N \mapsto {\bf 1}_M$ provides the require
isomorphism between $N$ and $[F, {\bf 1}_M]$. \hfill $\Box$

\bigskip

\noindent {\em Proof of Theorem \ref{t5.1}.} Let $M$ be reducible.
Let $A$ be a $\sim$-equivalence class of $M$ not contained in any
facet. Then $k = \#(A) \geq 2$.  Take $N$ as in Lemma \ref{l5.1}.
Then, clearly $M = N \ast S^{\,k-2}_k$, where $S^{\,k-2}_k$ is the
standard sphere with vertex set $A$. Conversely, if $M = N \ast
S^{\,k-2}_k(A)$ then any two vertices in $A$ are
$\sim$-equivalent. On the other hand, if $x\in A$ and $y \in V(N)$
then, by Lemma \ref{l2.1}, there is a facet $F$ of $N$ not
containing $y$. Hence $F\cup (A\setminus\{x\})$ is a facet of $M$
containing neither $x$ nor $y$. Therefore, by Lemma \ref{l4.2}
$(a)$, $x\not\sim y$. Thus $A$ is a $\sim$-equivalence class of
$M$. Since $A$ is not contained in any face of $M$, $M$ is
reducible.  \hfill $\Box$

\begin{defn}  \label{creducible}
{\rm A cellular pseudomanifold $M$ is said to be } completely
reducible {\rm if $M$ is the join of finitely many standard
spheres. (Thus, completely reducible cellular pseudomanifolds are
combinatorial spheres.)}
\end{defn}

From the proofs of Lemma \ref{l5.1} and Theorem \ref{t5.1} above,
it is immediate the following\,:

\begin{theo}  \label{t5.2}
Any cellular pseudomanifold is uniquely the join of an irreducible
cellular pseudomanifold and a completely reducible cellular pseudomanifold.
\end{theo}

\noindent {\em Proof.} Let $M$ be the given cellular
pseudomanifold. Suppose $M$ has $k$ $\sim$-equivalence classes
(say, $A_1, \dots, A_k$) which are not contained in any facet. Let
$S_i$ denote the standard spheres with vertex set $A_i$. Put $M_2
= S_1\ast \cdots \ast S_k$. Let $M_1$ denote the sub-poset of $M$
consisting of all faces of $M$ disjoint from $V(M_2)$ (with a
trivial face {\bf 1} thrown in). Then $M = M_1 \ast M_2$, $M_1$ is
irreducible and $M_2$ is completely reducible. Uniqueness of the
decomposition also follows from the proof of Lemma \ref{l5.1} (in
fact, the standard spheres in the completely reducible component
correspond to the $\sim$-equivalence classes $A_1, \dots, A_k$).
\hfill $\Box$

\begin{theo}  \label{t5.3}
The reducible cellular pseudomanifolds of excess $2$ and dimension
$d$ are $S^{\,b_1}_{b_1 + 2} \ast S^{\,b_2}_{b_2 + 2} \ast
S^{\,b_3}_{b_3 + 2} $ with $0 \leq b_1 \leq b_2 \leq b_3$, $b_1 +
b_2 + b_3 = d - 2$ and $((S^{\,d_1}_{d_1 + 2} \ast S^{\,d_2}_{d_2
+ 2}) \otimes S^{\,d_3}_{d_3 + 2}) \ast S^{\,d_4}_{d_4 + 2}$ with
$0 \leq d_1 \leq d_2$, $-1 \leq d_3$, $0 \leq d_4$, $d_1 + d_2 +
d_3 + d_4 = d - 4$. These are mutually non-isomorphic polytopal
spheres. Their number is $\lfloor \frac{(d^2 + 1)(2d - 1) + 9}{24}
\rfloor$.
\end{theo}

\noindent {\em Proof.} Let $M$ be a reducible $d$-dimensional
cellular pseudomanifold of excess 2. Then, by Theorem \ref{t5.1},
there exists a cellular pseudomanifold $N$ and an integer $c \geq
0$ such that $M = N \ast S^{\,c}_{c+2}$. By Lemma \ref{l2.5},
$\dim(N) = d-c-1$ and $e(N) = 1$. The first statement now follows
from Theorem \ref{t3.1}.

A cellular pseudomanifold of the first type is simplicial but a
cellular pseudomanifold of the second type is not simplicial.
Therefore, by Theorems \ref{t5.2} and \ref{t3.1}, all these
cellular pseudomanifolds are mutually non-isomorphic. It is known
(cf. \cite{bd2}) that all the cellular pseudomanifold of the first
type are polytopal spheres. By Theorem \ref{t3.1} and Theorem
\ref{t2.2} $(iii)$, all the cellular pseudomanifold of the second
type are polytopal spheres.

By Theorem \ref{t3.1}, the number of cellular pseudomanifolds of
the second type is
$$
\sum_{d_4=0}^{d-3} \left\lfloor\left(\frac{d - d_4 -
1}{2}\right)^2 \right\rfloor = \sum_{n=1}^{d-1} \left\lfloor
\frac{n^2}{4} \right\rfloor = \frac{1}{4} \sum_{n=1}^{d-1} n^2 -
\!\!\! \sum_{\begin{array}{c} {}^{n=1,} \\[-1mm] {}^{n ~ {\rm odd}}
\end{array}}^{d-1}\!\! \frac{1}{4} = \frac{(d-1)d(2d-1)}{24}
- \frac{1}{4}\left\lfloor\frac{d}{2}\right\rfloor.
$$
From the classification of excess 2 pseudomanifolds in
\cite{bd2}, we know that the number of cellular pseudomanifolds of
first type is $\left\lfloor\frac{(d+4)(d-2)}{12}\right\rfloor +1$.
Therefore, the total number is $\left \lfloor \frac{(d + 4)(d -
2)}{12} \right \rfloor + 1 + \frac{(d - 1)d(2d - 1)}{24} -
\frac{1}{4} \left \lfloor \frac{d}{2} \right \rfloor = \left
\lfloor \frac{2d^3 - d^2 + 2d +8}{24} \right \rfloor = \left
\lfloor \frac{(d^2 +1)(2d-1) + 9}{24} \right \rfloor$. \hfill
$\Box$

\begin{cor}  \label{c5.1}
The reducible neighbourly cellular pseudomanifolds of excess $2$
and dimension $d$ are $S^{\,b_1}_{b_1 + 2} \ast S^{\,b_2}_{b_2 +
2} \ast S^{\,b_3}_{b_3 + 2}$ with $1 \leq b_1 \leq b_2 \leq b_3$,
$b_1 + b_2 + b_3 = d - 2$ and $((S^{\,d_1}_{d_1 + 2} \ast
S^{\,d_2}_{d_2 + 2}) \otimes S^{\,d_3}_{d_3 + 2}) \ast
S^{\,d_4}_{d_4 + 2}$ with $1 \leq d_1 \leq d_2$, $-1 \leq d_3$, $1
\leq d_4$, $d_1 + d_2 + d_3 + d_4 = d - 4$. These are mutually
non-isomorphic polytopal spheres. Their number is $\left \lfloor
\frac{(d^2 - 6d - 8)(2d - 7) + 9}{24} \right \rfloor$.
\end{cor}

\noindent {\em Proof.} Observe that $S^{\,b_1}_{b_1 + 2} \ast
S^{\,b_2}_{b_2 + 2} \ast S^{\,b_2}_{b_2 + 2}$ is neighbourly if
and only if $b_1, b_2, b_3 \geq 1$ and $((S^{\,d_1}_{d_1 + 2} \ast
S^{\,d_2}_{d_2 + 2}) \otimes S^{\,d_3}_{d_3 + 2}) \ast
S^{\,d_4}_{d_4 + 2}$ is neighbourly if and only if $d_1, d_2, d_4
\geq 1$. Thus, the first and second parts follow from Theorem
\ref{t5.3}.

So, the number of reducible neighbourly cellular pseudomanifolds
of excess 2 and dimension $d$ is $= \#(\{(b_1, b_2, b_3) ~ : ~ 1
\leq b_1 \leq b_2 \leq b_3$, $b_1 + b_2 + b_3 = d - 2\})+
\#(\{(d_1, d_2, d_3, d_4) ~ : 1 \leq d_1 \leq d_2, ~ -1 \leq d_3,
~ 1 \leq d_4, ~ d_1 + d_2 + d_3 + d_4 = d - 4\}) = \#(\{(a_1, a_2,
a_3) ~ : ~ 0 \leq a_1 \leq a_2 \leq a_3$, $a_1 + a_2 + a_3 = d -
5\})+ \#(\{(c_1, c_2, d_3, c_4) ~ : 0 \leq c_1 \leq c_2, ~ -1 \leq
d_3, ~ 1 \leq c_4, ~ c_1 + c_2 + d_3 + c_4 = d - 7\})$. By the
proof of Theorem \ref{t5.3}, this number is $= \left \lfloor
\frac{((d-3)^2 +1)(2(d-3)-1) + 9}{24} \right \rfloor$. This
completes the proof. \hfill $\Box$

\section{\sc Proper cellular pseudomanifolds}

\begin{defn}  {\rm A cellular pseudomanifold $M$ is
called {\em proper} if each $\sim$-equivalence class of $M$ is
(the shadow of) a face of $M$.  We shall say that $M$ is {\em
improper} otherwise. }
\end{defn}

\begin{defn}  {\rm A cellular pseudomanifold is
called {\em primitive} if for any two vertices $x$, $y$ of $M$,
$x\sim y$ implies $x = y$, i.e., if all the $\sim$-equivalence
classes of $M$ are singletons. Thus any primitive cellular
pseudomanifold is proper.}
\end{defn}

$S^{\,1}_5 \otimes S^{\,- 1}_1$, $S_2$, $S_3$, $S_4$, $R$ and
$\mathbb{R}\mbox{P}^{\,2}_6$ in Example \ref{v6d2} are primitive. $S_1$
is proper. All proper cellular pseudomanifolds are irreducible,
but the converse is false. For instance, $(S^{\,0}_2 \ast
S^{\,0}_2) \otimes S^{\,-1}_1$ and $S_5$ in Example \ref{v6d2} are
irreducible but not proper. However, `irreducible' and `proper'
become equivalent notion within the class of normal 
pseudomanifolds.

\begin{nota}  {\rm Let $M$ be a proper cellular
pseudomanifold.  Let $M/\!\!\sim$ denote the poset whose elements
are the sets $\{X_1, \dots, X_t\}$, where $X_1, \dots, X_t$ are
$\sim$-equivalence classes of vertices of $M$ such that $X_1 \cup
\cdots \cup X_t$ is (the shadow of) a face of $M$.  The partial
order on $M/\!\!\sim$ is set inclusion.}
\end{nota}

\begin{eg}  \label{nonproperex} {\rm Let $M=
(S^{\,d_1}_{d_1 +2}(V_1) \ast S^{\,d_2}_{d_2 +2}(V_2)) \otimes
S^{\,d_3}_{d_3 +2}(V_3)$ for some $d_1, d_2, d_3$ with $0\leq d_1
\leq d_2$, $-1 \leq d_3$, $d_1 + d_2 + d_3 = d-3$ as in Theorem
\ref{t3.1}. Then $V_1$, $V_2$, $V_3$ are the $\sim$-equivalence
classes in $M$ and $V_1$, $V_2$ are not faces but $V_1\cup V_2$,
$V_3$ are faces.}
\end{eg}

\begin{theo}  \label{t6.1} For any proper cellular
pseudomanifold $M$, $M/\!\!\sim$ is a primitive cellular
pseudomanifold with $e(M/\!\!\sim) = e(M)$.
\end{theo}

\noindent {\em Proof.} Let $A_1, \dots, A_m$ be the
$\sim$-equivalence classes and let $f$ be as in Lemma \ref{l4.3}.
For $U\subseteq V(M)$, let $\widetilde{U}$ denote the collection
of all the $\sim$-equivalence classes which are contained in $U$
and $\widehat{U}$ is as defined in Lemma \ref{l4.2}. So,
$\widehat{U}$ is the union of all the members of $\widetilde{U}$.

Since $M$ is proper, each $A_i$ is a face and hence $V(M/\!\!\sim)
= \{A_1, \dots, A_m\}$. Again, $A_i$ is a face implies, by part
$(c)$ of Lemma \ref{l4.2}, $A_i\cup f$ is a face for each $i$ and
hence $V({\rm lk}(f)) = \{f\cup A_1, \dots, f\cup A_m\}$. Now, a
collection ${\cal A}$ of $\sim$-equivalence classes of $M$ forms a
face of $M/\!\!\sim$ if and only if $\cup{\cal A}$ is a face of
$M$. By part $(c)$ of Lemma \ref{l4.2}, this holds if and only if
$(\cup {\cal A}) \cup f$ is a face of $M$. Thus, the map $\psi
\colon M/\!\!\sim \, \to {\rm lk}(f)$ defined by $\psi({\cal A}) =
(\cup {\cal A}) \cup f$ is an isomorphism and hence $M/\!\!\sim$
is a cellular pseudomanifold.

If possible, let $A_i\sim A_j$ in $M/\!\!\sim$ for $i\neq j$.
Then, by part $(a)$ of Lemma \ref{l4.2}, each facet in
$M/\!\!\sim$ has $A_i$ or $A_j$ as a vertex. This implies that
each facet in $M$ has $A_i$ or $A_j$ as a subset. Now,
$x_i\not\sim x_j$ in $M$ implies that there exists a facet
$\alpha$ in $M$ such that $x_i\in \alpha$, $x_j\not\in \alpha$ and
$\beta = (x_i, x_j)(\alpha)$ is not a facet in $M$. Then, by part
$(c)$ of Lemma \ref{l4.2}, $\widehat{\alpha}$ is a face of $M$ but
$\widehat{\beta}$ is not a face of $M$. This implies, by the
definition of $M/\!\!\sim$, that $\widetilde{\alpha}$ is a face of
$M/\!\!\sim$ and $\widetilde{\beta}$ is not a face of
$M/\!\!\sim$.

Since $x_j\not\in\alpha$, $A_j\not\subseteq \alpha$ and hence
$A_i\subseteq \alpha$. Since $\alpha$ is a facet and $x_j\not\in
\alpha$, $A_j\setminus\{x_j\}\subseteq \alpha$ and hence $\beta =
(\alpha \setminus A_i) \cup (A_i\setminus \{x_i\}) \cup A_j$. This
implies $\widehat{\beta}= (\widehat{\alpha}\setminus A_i) \cup
A_j$ and hence $\widetilde{\beta} = (\widetilde{\alpha} \setminus
\{A_i\}) \cup \{A_j\} = (A_i, A_j)(\widetilde{\alpha})$. This is a
contradiction since $\widetilde{\alpha}$ is a face of
$M/\!\!\sim$, $\widetilde{\beta}$ is a not face of $M/\!\!\sim$
and $(A_i, A_j)$ is an automorphism of $M/\!\!\sim$. So,
$A_i\not\sim A_j$. This implies that $M/\!\!\sim$ is primitive.

Let $M$ have $n$ vertices. Then $\#(f) = n - m$ and hence, by part
$(a)$ of Lemma \ref{l4.1}, $f$ is a face of dimension $n-m-1$.
This implies $\dim(M/\!\!\sim) = \dim({\rm lk}(f)) = d+1 - (n-m-1)
- 2 = d-n+m$ and hence $e(M/\!\!\sim) = m - (d-n+m) -2 = n-d-2 =
e(M)$. \hfill $\Box$

\begin{nota}  {\rm Let $N$ be a proper cellular
pseudomanifold and $f\colon V(N)\to \mathbb{Z}^{+}$ be a map. Let $\{V_x
\, \colon \, x\in V(N)\}$ be a family of pairwise disjoint sets,
indexed by the vertices of $N$, such that $\#(V_x) = f(x)$ for all
$x\in V(N)$. Let $\langle N, f\rangle$ denote the poset (with
set-inclusion as the partial order) whose members are the subsets
$A$ of $\cup_{x\in V(N)} V_x$ such that $\{x\in V(N) \, \colon \,
V_x\subseteq A\}$ is the shadow of a face of $N$. }
\end{nota}

\begin{lemma}  \label{l6.1}
Let $N$ be a proper cellular pseudomanifold and let $f\colon V(N)
\to \mathbb{Z}^{+}$ be a map. Then $M := \langle N, f\rangle$ is a proper
cellular pseudomanifold with $V(M) = \cup_{x\in V(N)} V_x$. For
each $\sim$-equivalence class $A$ of $N$, $\bar{A} := \cup_{x\in
A} V_x$ is a $\sim$-equivalence class in $M$. If further, $N$ is
primitive then $M/\!\!\sim$ is isomorphic to $N$.
\end{lemma}

\noindent {\em Proof. } As usual, we identify a face of $N$ with
its shadow. Put $V = \sqcup_{x\in V(N)} V_x$. For any subset $A$
of $V$, let $\widetilde{A} = \{x\in V(N) ~ \colon ~ V_x\subseteq
A\}$. Then $A$ is a face of $M = \langle N, f\rangle$ if and only
if $\tilde{A}$ is a face. Define $\rho_{{}_M}\colon M\to \mathbb{N}$ by
$$
\rho_{{}_M}(A) = \#(A) - \#(\widetilde{A}) +
\rho_{{}_N}(\widetilde{A})
$$
where $\rho_{{}_N}$ is the rank function of $N$. Then it is easily
verified that $M$ is a ranked poset with rank function
$\rho_{{}_M}$.

Let $A$, $B$ be two faces of $M$ with $A\subseteq B$. Then
$\widetilde{A}\subseteq \widetilde{B}$ and, from the formula of
$\rho_{{}_M}$, it follows that
\begin{eqnarray} \label{rankeqn}
\rho_{{}_M}(B) - \rho_{{}_M}(A) =
\rho_{{}_N}(\widetilde{B}) - \rho_{{}_N}(\widetilde{A}) +
\!\!\sum_{x\in\widetilde{B}\setminus \widetilde{A}}
(\#(V_x\setminus A) - 1) +
\sum_{x\not\in\widetilde{B}}\#(V_x\cap(B\setminus A)).
\end{eqnarray}

Now, if $A$, $B$ are two faces of $M$ with $A\subset B$ and
$\rho_{{}_M}(B) - \rho_{{}_M}(A) = 2$ then $\widetilde{A}
\subseteq \widetilde{B}$ and $\rho_{{}_N}(\widetilde{B}) -
\rho_{{}_M}(\widetilde{A}) \leq 2$. First suppose that
$\rho_{{}_N}(\widetilde{B}) - \rho_{{}_M}(\widetilde{A}) = 2$.
Then (\ref{rankeqn}) implies that for $x\not\in \widetilde{B}$,
$V_x$ is disjoint from $B\setminus A$ and for $x \in \widetilde{B}
\setminus \widetilde{A}$, $A$ misses exactly one element of $V_x$.
Therefore, if $\widetilde{C}_i$, $i= 1, 2$ are the faces of $N$
strictly between $\widetilde{A}$ and $\widetilde{B}$ then $C_i :=
A\cup(\cup_{x\in\widetilde{C}_i} V_x)$, $i=1, 2$ are the only
faces of $M$ strictly between $A$ and $B$.

Next suppose that $\rho_{{}_N}(\widetilde{B}) -
\rho_{{}_M}(\widetilde{A}) = 1$. Then no face of $N$ lies strictly
between $\widetilde{A}$ and $\widetilde{B}$. Therefore, for any
face $C$ of $M$ strictly between $A$ and $B$, either
$\widetilde{C} = \widetilde{A}$ or $\widetilde{C} =
\widetilde{B}$. Now, (\ref{rankeqn}) shows that there are two
cases. In Case 1, there is a unique $x_0\in \widetilde{B}
\setminus \widetilde{A}$ such that $B\setminus A = V_{x_0}
\setminus A$ is a doubleton, say $B\setminus A = \{u, v\}$. In
this case, $A\cup\{u\}$, $A\cup\{v\}$ are the only two faces of
$M$ strictly between $A$ and $B$. In Case 2, there is a unique
$x_0\not\in\widetilde{B}$ such that $V_{x_0}\cap (B\setminus A)$
is a singleton, say $V_{x_0}\cap (B\setminus A) = \{w\}$. Then
$C_1 := A\cup (\cup_{x\in \widetilde{B}} V_x)$ and $C_2 := A \cup
\{w\}$ are the only two faces of $M$ strictly between $A$ and $B$.
(Observe that, $\widetilde{C}_1 = \widetilde{B}$ and
$\widetilde{C}_2 = \widetilde{A}$ in the last case.)

Finally, suppose $\rho_{{}_N}(\widetilde{B}) -
\rho_{{}_N}(\widetilde{A}) = 0$, i.e., $\widetilde{B} =
\widetilde{A}$. By (\ref{rankeqn}), $B\setminus A$ is a doubleton,
say $B\setminus A = \{u, v\}$. In this case, $A\cup\{u\}$ and $A
\cup \{v\}$ are the only faces of $M$ strictly between $A$ and
$B$. So, in all cases, if $A$, $B$ are faces of $M$ with $A\subset
B$ and $\rho_{{}_M}(B) - \rho_{{}_M}(A) = 2$ then there are
exactly two faces of $M$ strictly between $A$ and $B$.

Now, let $A$, $B$ be faces of $M$ with $A \subseteq B$ and
$\rho_{{}_M}(B) - \rho_{{}_M}(A) \geq 3$. First suppose that
$\rho_{{}_N}(\widetilde{B}) - \rho_{{}_N}(\widetilde{A}) =
\rho_{{}_M}(B) - \rho_{{}_M}(A)$. Then for $A < C < B$,
$\rho_{{}_N}(\widetilde{B}) - \rho_{{}_N}(\widetilde{C}) \leq
\rho_{{}_M}(B) - \rho_{{}_M}(C)$ and $\rho_{{}_N}(\widetilde{C}) -
\rho_{{}_N}(\widetilde{A}) \leq \rho_{{}_M}(C) - \rho_{{}_M}(A)$.
This implies that
\begin{eqnarray} \label{=eqn}
\rho_{{}_N}(\widetilde{B}) - \rho_{{}_N}(\widetilde{C}) =
\rho_{{}_M}(B) - \rho_{{}_M}(C) & \mbox{and} &
\rho_{{}_N}(\widetilde{C}) - \rho_{{}_N}(\widetilde{A}) =
\rho_{{}_M}(C) - \rho_{{}_M}(A).
\end{eqnarray}
For $\alpha\in [\widetilde{A}, \widetilde{B}]$, if $C :=
A\cup(\cup_{x\in\alpha} V_x)$ then $C\in [A, B]$ and
$\widetilde{C}=\alpha$. Thus, $C\mapsto \widetilde{C}$ is a map
from $[A, B]$ onto $[\widetilde{A}, \widetilde{B}]$. Also, because
of (\ref{=eqn}), this map is one-to-one. Thus, it defines an
isomorphism between $[A, B]$ and $[\widetilde{A}, \widetilde{B}]$.
Hence $[A, B]$ is a cellular pseudomanifold. In particular,
$\La([A, B])$ is connected.

Next suppose that $\rho_{{}_N}(\widetilde{B}) -
\rho_{{}_N}(\widetilde{A}) < \rho_{{}_M}(B) - \rho_{{}_M}(A)$. If
$B$ is a union of $\sim$-equivalence classes then for any facet
$C$ of $[A, B]$,  $\widetilde{C}\neq \widetilde{B}$ and hence
$\widetilde{C}$ is a facet of $[\widetilde{A}, \widetilde{B}]$.
Fix $\bar{x}\in V_x\setminus A$ for $x\in \widetilde{B}\setminus
\widetilde{A}$. Define $\varphi\colon [A, B] \to [A, B]$ by
$\varphi(D) = (\cup_{x\in \widetilde{D}} V_x)\cup (\cup_{x\in
\widetilde{B} \setminus \widetilde{D}} V_x \setminus\{\bar{x}\})$.
Now, if $C$ is a facet of $[A, B]$ then $C$ misses exactly one
vertex from $V_x$ for $x \in \widetilde{B} \setminus
\widetilde{C}$ and either $C = \varphi(C)$ or $C$ is adjacent to
$\varphi(C)$ in $\La([A, B])$. If $\alpha$ is a face of
$[\widetilde{A}, \widetilde{B}]$ then define $\widehat{\alpha} =
(\cup_{x \in \alpha} V_x) \cup (\cup_{x \in \widetilde{B}
\setminus \alpha} V_x\setminus \{\bar{x}\})$. Clearly, $\alpha
\mapsto \widehat{\alpha}$ defines an isomorphism from
$[\widetilde{A}, \widetilde{B}]$ to $\varphi([A, B])$. Now, let
$C$ and $D$ be two facets in $[A, B]$. Let $\widetilde{C} =
\alpha_1 \cdots \alpha_k =\widehat{D}$ be a path in
$\La([\widetilde{A}, \widetilde{B}])$. ($[\widetilde{A},
\widetilde{B}]$ is a normal pseudomanifold and hence $\La([\widetilde{A},
\widetilde{B}])$ is connected.) Let $P$ be a path from $C$ to
$\varphi(C)=\widehat{\alpha}_1$ and $Q$ be a path from
$\widehat{\alpha}_k = \varphi(D)$ to $D$ in $\La([A, B])$. (If
$C=\varphi(C)$ then $P$ consists of a single vertex.) Then $P
\widehat{\alpha}_2\cdots \widehat{\alpha}_{k-1}Q$ is a path from
$C$ to $D$ in $\La([A, B])$. Thus, $\La([A, B])$ is connected.

Thus, assume that $B$ is not union of $\sim$-equivalence classes.
Then there exists $x_0\in V(N)\setminus\widetilde{B}$ such that
$V_{x_0}\cap B\neq\emptyset$. If $C$ is a facet in $[A, B]$ then
either $\widetilde{C}$ is a facet of $[\widetilde{A},
\widetilde{B}]$ or $\widetilde{C} = \widetilde{B}$. Let ${\cal C}
:= \{C ~ \colon ~ C$ a  facet of $[A, B]$ and $\widetilde{C} =
\widetilde{B}\}$. If $C \in {\cal C}$ then, by (\ref{rankeqn}),
there exists unique $x_1\not\in \widetilde{B} = \widetilde{C}$
such that $V_{x_1}\cap (B\setminus C)$ is a singleton, say,
$V_{x_1}\cap (B\setminus C)= \{w\}$. So, $B=C\cup\{w\}$ or
$C=B\setminus\{w\}$, where $w\in V_{x_1}\cap B \neq V_{x_1}$. This
implies that any two elements of $\cal C$ are adjacent in $\La([A,
B])$. If $D$ is a facet of $[A, B]$ such that $D\not\in {\cal C}$,
i.e., $\widetilde{D}$ is a facet of $[\widetilde{A},
\widetilde{B}]$ then, by (\ref{rankeqn}), $V_{x}\cap (B\setminus
D) = \emptyset$ for any $x\not\in\widetilde{B}$. So, $V_x\cap B =
V_x\cap D$ for any $x\not\in\widetilde{B}$. In particular,
$V_{x_0} \cap D = V_{x_0}\cap B\neq\emptyset$. Let $a\in V_{x_0}
\cap D = V_{x_0}\cap B$. Let $C=B\setminus\{a\}$. Then $C$ is a
face of $M$ and hence a facet of $[A, B]$. Clearly,
$\widetilde{C}= \widetilde{B}$ and hence $C\in {\cal C}$. Since
$\#(C\cap D) = \#(D) -1$, $C$ and $D$ are adjacent in $\La([A,
B])$. Thus each vertex of $\La([A, B])$ outside $\cal C$ is
adjacent to a vertex inside $\cal C$ and any two vertices inside
$\cal C$ are adjacent. This implies that $\La([A, B])$ is
connected. This proves that $M$ is a cellular pseudomanifold.

Let $x\not\sim y$ in $N$ and $u\in V_x$, $v\in V_y$. Then there is
a facet $\tau$ of $N$ such that $x\in \tau$, $y\not\in\tau$ (say)
and $(\tau\setminus\{x\})\cup\{y\}$ is not facet of $N$. Put
$\alpha = \cup_{z\in\tau\setminus\{x\}} V_z$, $\beta = \alpha \cup
(V_x\setminus\{u\}) \cup (V_y\setminus\{v\})$. Then $\beta\cup
\{u\}$ is a face of $M$ but $\beta\cup \{v\}$ is not a face of
$M$. Therefore $u\not\sim v$ in $M$. So, for $x\not\sim y$ in $N$,
an element of $V_x$ is not equivalent to an element of $V_y$ in
$M$. Now, let $A$ be a $\sim$-equivalence class of vertices of
$N$. Put $\bar{A} = \cup_{x\in A} V_x$. From above argument, it
follows that if $u\in A$ and $v\not\in A$ then $u\not\sim v$ in
$M$. Again, from the definition of $M$, it is known that the
interchange of any two vertices in $\bar{A}$ is an automorphism of
$M$. That is, any two vertices in $\bar{A}$ are equivalent in $M$.
Therefore, $\bar{A}$ is an $\sim$-equivalence class of vertices of
$M$. This establishes the validity of the definition of the
$\sim$-equivalence classes of $M$. From this description, it is
immediate that all the $\sim$-equivalence classes of $M$ are
faces, so that $M$ is proper.

If $N$ is primitive then it follows that the $\sim$-equivalence
classes of $M$ are $V_x$, $x\in V(N)$. That is, the vertices of
$M/\!\!\sim$ are $V_x$, $x\in V(N)$. From the definition of
$M/\!\!\sim$, it is now immediate that $x \mapsto V_x$ is the
require isomorphism between $N$ and $M/\!\!\sim$ in this case.
\hfill $\Box$

\begin{lemma}  \label{l6.2}
For $i= 1, 2$, let $N_i$ be primitive cellular pseudomanifolds and
$f_i \colon V(N_i) \to \mathbb{Z}^{+}$. Then $\langle N_1, f_1 \rangle$
and  $\langle N_2, f_2 \rangle$ are isomorphic if and only if
there is an isomorphism $\varphi \colon N_2 \to N_1$ such that
$f_2 = f_1 \circ \varphi$.
\end{lemma}

\noindent {\em Proof.} Suppose $\varphi \colon N_1 \to N_2$ be an
isomorphism such that $f_2(x) = f_1(\varphi(x))$ for $x \in
V(N_1)$. Let $V_x$, $x \in V(N_1)$ be the $\sim$-equivalence
classes of $M_1 = \langle N_1, f_1 \rangle$ and $W_y$, $y \in
V(N_2)$ be the $\sim$-equivalence classes of $M_2 = \langle N_2,
f_2 \rangle$. Thus $f_1(x) = \#(V_x)$, $f_2(y) = \#(W_y)$. Since
$f_2 = f_1 \circ \varphi$, there is a bijection $h_x \colon V_x
\to W_{\varphi(x)}$, for $x\in V(N_1)$. Then $\cup_{x\in V(N_1)}
h_x$ is the require isomorphism between $M_1$ and $M_2$.
Conversely, if $M_1$ and $M_2$ are isomorphic, then any
isomorphism between them induces a bijection $\varphi$ between the
set of $\sim$-equivalence classes of $M_1$ and the set of
$\sim$-equivalence classes of $M_2$. Then $\varphi \colon V(N_1)
\to V(N_2)$ (under the identification given in Lemma \ref{l6.1})
is the require isomorphism between $N_1$ and $N_2$ intertwining
$f_1$ and $f_2$. \hfill $\Box$

\begin{theo}  \label{t6.2}
Let $M$ be a proper cellular pseudomanifold. Then $M$ is
isomorphic to $\langle N, f\rangle$, where $N=M/\!\!\sim$ is
primitive and $f\colon V(N) \to \mathbb{Z}^+$ is given by $f(A) = \#(A)$
for each $\sim$-equivalence class $A$ of $M$. Further, this
representation of $M$ is unique up to isomorphism.
\end{theo}

\noindent {\em Proof.} The result follows from Lemmas \ref{l6.1}
and \ref{l6.2}. \hfill $\Box$

\section{\sc  Cellular spheres of excess 2}

As before, we identify a face of a cellular pseudomanifold with
its shadow.

We briefly recall the theory of  Gale diagram here. For further
details, see \cite{g}. If $P$ is an $n$-vertex polytopal sphere of
excess $e>0$, then a {\em Gale diagram} of $P$ is a multi-set $G$
of $n$ points in $S^{\,e-1} \cup \{0\}$ (where as usual, $S^{\,e
- 1}$ is the unit sphere in $\mathbb{R}^{e}$, $0$ is the center of
$S^{\,e - 1}$) together with a surjection $g \colon V(P) \to G$
such that $A \subseteq V(P)$ is a co-facet ($\equiv$ complement
of a facet) if and only if $\#(g(A)) = \#(A) \leq e+1$ and $0$
belongs to the relative interior of the convex hull of $f(A)$.
Every polytopal sphere has a Gale diagram. Conversely, any
multi-set $G$ of points in $S^{\,e-1}\cup\{0\}$ is a Gale diagram
of a polytopal sphere of excess $e$ provided, for each hyperplane
$\pi$ in $\mathbb{R}^{e}$ passing through $0$ and a point of $G$, each of
the open hemisphere determined by $\pi$ contains at least two
points of $G$. Two Gale diagrams are called {\em isomorphic} if
the corresponding polytopal spheres are isomorphic.

\begin{lemma}  \label{l8.1}
Let $P$ be a polytopal $d$-sphere of excess $2$ with a Gale
diagram $G$ $(\subseteq S^{\,1}\cup \{0\})$. Let $\alpha$ and
$\beta$ be two distinct points in $G$ corresponding to the
vertices $x$ and $y$ of $P$. Let $H$ be the Gale diagram obtained
from $G$ by shifting $\alpha$ to $\beta$. \vspace{-1mm}
\begin{enumerate}
\item[$(i)$] If $x \sim y$ in $P$ then $G$ has no point in the
shorter closed arc $[-\alpha, -\beta]$. \vspace{-1mm}
\item[$(ii)$] If $\alpha = \beta$ then $x\sim y$ in $P$.
\vspace{-1mm}
\item[$(iii)$] If $[-\alpha, -\beta]$ does not contain any point
of $G$ then $G$ and $H$ are isomorphic.
\end{enumerate}
\end{lemma}

\noindent {\em Proof.} Let $x \sim y$ in $P$. If possible let
there exist a point $\gamma$ of $G$ in $[-\alpha, -\beta]$. Let
$z$ be a vertex corresponding to $\gamma$. If $\gamma = -\alpha$
then $0$ is in the relative interior of the convex hull of
$\{\alpha, \gamma\}$ (the line segment joining $\alpha$ and
$\gamma$) but not in convex hull of $\{\beta, \gamma\}$. So, $V(P)
\setminus \{z, x\}$ is a facet of $P$ but $V(P) \setminus \{z,
y\}$ is not a facet of $P$. This is not possible, since $x\sim
y$. So, $\gamma \neq - \alpha$. Similarly, $\gamma \neq - \beta$.
So, assume that $\gamma \in (- \alpha, - \beta)$. Then the
interior of the convex hull of $\{\alpha, \beta, \gamma\}$
contains $0$ and hence $V \setminus \{x, y, z\}$ is a facet. This
is not possible by Lemma \ref{l4.2} $(a)$. This proves $(i)$.

Let $\alpha = \beta$. Let $\sigma$ be a facet containing $x$ but
not containing $y$. This implies that there exists a subset $X$
of $G$ with at most three points such that, $\beta\in X$, $\alpha
\not\in X$ and  the relative interior of the convex hull of $X$
contains $0$ ($X$ correspond to the complement of $\sigma$).
Clearly, the convex hull of $X$ is same as the convex hull of $(X
\setminus \{\beta\}) \cup \{\alpha\}$. Therefore, $(\sigma \cup
\{y\}) \setminus \{x\}$ is a facet of $P$. This implies that $x
\sim y$ in $P$. This proves $(ii)$.

$(iii)$ follows from Theorem 1 in \cite[Section 6.3]{g}. \hfill
$\Box$

\begin{lemma}  \label{l8.2}
Let $M$ be a $d$-dimensional cellular pseudomanifold on $d+4$
vertices such that $M$ has no facet containing $d+3$ vertices.
\begin{enumerate}
\item[$(i)$] If $\alpha$ is a $c$-face with $c+2$ vertices then
the link of $\alpha$ is of the form $S^{\,r}_{r + 2} \ast
S^{\,s}_{s + 2}$ for some $r \geq s \geq - 1$ with $r + s = d - c
-2$.
\item[$(ii)$] If $M$ has no face $\tau$ such that $\partial_M \tau =
S^{\,a}_{a + 2} \ast S^{\,b}_{b + 2}$ for some $a, b \geq 0$ then
$M$ is a normal pseudomanifold.
\end{enumerate}
\end{lemma}

\noindent {\bf Proof.} Since ${\rm lk}_M(\alpha)$ is
$(d-c-1)$-dimensional and there are $d-c+2$ vertices outside
$\alpha$, the excess of ${\rm lk}_M(\alpha)$ is 0 or 1. If the
excess is 0 then ${\rm lk}_M(\alpha) = S^{\,d-c-1}_{d-c+1} =
S^{\,d-c-1}_{d-c+1} \ast S^{\,-1}_{1}$. So, assume that the
excess of ${\rm lk}_M(\alpha)$ is 1. If possible let
$(S^{\,d_1}_{d_1 + 2}(D_1) \ast S^{\,d_2}_{d_2 + 2}(D_2)) \otimes
S^{\,d_3}_{d_3 + 2}(D_3)$ for some $d_2 \geq d_1 \geq 0$, $d_3
\geq - 1$ with $d_1 + d_2 + d_3 = d - c - 4$. Then, for $z\in
D_3$, $\alpha \cup D_1 \cup D_3 \cup (D_3 \setminus \{z\})$ is a
facet of $M$ with $d+3$ vertices. A contradiction. Therefore, by
Theorem \ref{t3.1}, ${\rm lk}_M(\alpha)$ is a join of two
standard spheres. This proves $(i)$.

Since $M$ has no facet containing $d+3$ vertices, $M$ has no
$i$-face with $i+3$ vertices for all $i \leq d$. If possible let
$M$ has an $c$-face $\tau$ with $c + 2$ vertices for some $c \leq
d$. Since $M$ has no face whose boundary is the join of two
spheres, by Theorem \ref{t3.1}, $\partial_{M} \tau =
(S^{\,d_1}_{d_1 +2}(D_1) \ast S^{\,d_2}_{d_2 +2}(D_2)) \otimes
S^{\,d_3}_{d_3 +2}(D_3)$ for some $d_1$, $d_2$, $d_3$ with $d_2
\geq d_1 \geq 0$, $d_3 \geq - 1$ and $d_1 + d_2 + d_3 = c - 4$.
Then $\alpha := D_1 \cup D_2$ is a face of $M$ with boundary
$S^{\,d_1}_{d_1 + 2}(D_1) \ast S^{\,d_2}_{d_2 + 2}(D_2)$. This is
not possible. So, $M$ has no $i$-face with more than $i + 1$
vertices for all $i$. This implies that $M$ is a simplicial
complex. This proves $(ii)$ \hfill $\Box$

\begin{lemma}  \label{l8.3}
Let $M$ be a $d$-dimensional cellular pseudomanifold on $d+4$
vertices such that $M$ has no facet containing $d+3$ vertices.
Let $\tau$ be a $k$-face with $k+2$ vertices and $\partial \tau =
S^{\,a}_{a + 2}(A) \ast S^{\,b}_{b + 2}(B)$ for some $a, b \geq
0$. Let ${\rm lk}_M(\tau) = S^{\,r}_{r + 2}(R) \ast S^{\,s}_{s +
2}(S)$ for some $r \geq s \geq - 1$ with $r + s = d - k -2$. Let
$M_1$ be the cellular pseudomanifold whose facet-set is
$$
\{\sigma : \sigma \mbox{ facet in } M, \tau\not\subseteq \sigma\}
\cup \{A \cup (B \setminus \{y\}) \cup (R\setminus \{z\}) \cup
(S\setminus\{w\}) : y \in B, z \in R, w \in S\}.
$$
\begin{enumerate}
\item[$(i)$] $B(M_1)$ and $B(M)$ are combinatorially equivalent. \vspace{-1mm}
\item[$(ii)$] $M$ is a cellular sphere if and only if $M_1$ is so.
\vspace{-1mm}
\item[$(iii)$] If $M_1$ is a polytopal sphere then $M$ is so.
\end{enumerate}
\end{lemma}

\begin{lemma}  \label{l8.4}
Let all the notations be as in Lemma $\ref{l8.3}$. If $M_1$ is a
polytopal sphere then $M_1$ has a Gale diagram $G_1$ in which all
the vertices in $R$ are represented by one point $\rho \in G_1$,
all the vertices in $S$ are represented by one point $\sigma \in
G_1 \setminus \{\rho, -\rho\}$, all the vertices in $A$ are
represented by points of $G_1$ in the shorter open arc $(\rho,
\sigma)$ joining $\rho$ and $\sigma$, and all the vertices in $B$
are represented by points of $G_1$ in the shorter open arc
$(-\rho, -\sigma)$ joining $- \rho$ and $- \sigma$.
\end{lemma}

\noindent {\em Proof.} Observe that there is no face of $M_1$
containing $\tau$.

\smallskip

\noindent {\em Claim} 1. $B$ is not contained in any face of $M_1$
and ${\rm lk}_{M_1}(A) = S^{\,b}_{b + 2}(B) \ast S^{\,r}_{r +
2}(R) \ast S^{\,s}_{s + 2}(S)$.

\smallskip

Let $\alpha$ be a face of $M$ containing $B$. If $\tau \not
\subseteq \alpha$ then $\beta := \alpha \cap \tau$ is a proper
face of $\tau$ containing $B$ in $M$. Hence $B \subseteq \beta \in
\partial_M \tau$, a contradiction to the fact that $\partial_M
\tau = S^{\,a}_{a + 2}(A) \ast S^{\,b}_{b + 2}(B)$). So, $\tau
\subseteq \alpha$. Thus, $\tau$ is the smallest face of $M$
containing $B$. This implies that $B$ is not contained in any
facet of $M_1$. This proves the first part.

By the same argument, $\tau$ is the smallest face of $M$
containing $A$. This implies that ${\rm lk}_{M_1}(A) = S^{\,b}_{b
+ 2}(B) \ast S^{\,r}_{r + 2}(R) \ast S^{\,s}_{s + 2}(S)$. This
proves the claim.

\smallskip

\noindent {\em Claim} 2. All the vertices in $R$ are mutually
$\sim$-equivalent in $M_1$. (And similarly, all the vertices in
$S$ are mutually $\sim$-equivalent in $M_1$.)

\smallskip

To prove this claim, it suffices to show that if $z, z^{\,\prime}
\in R$ are distinct vertices and $C$ is a facet of $M_1$ such
that $z^{\,\prime} \in C$, $z \not\in C$, then $(C \cup \{z\})
\setminus \{z^{\,\prime}\}$ is also a facet of $M_1$.

If $A \subseteq C$ then from the known structure of ${\rm
lk}_{M_1}(A)$ it is clear that $C \cup \{z\} \setminus
\{z^{\,\prime}\}$ is also a facet of $M_1$. So, assume that $A
\not \subseteq C$. Let $x \in A \setminus C$.

By Claim 1, $B \not \subseteq C$. So, there exists $y \in B
\setminus C$. Since $M_1$ is of excess 2, it follows that $C =
V(M_1) \setminus \{x, y, z\}$ is a simplex. Fix an element $w \in
S$. Since $C_1 := C \setminus \{w\} = V(M_1) \setminus \{x, y, z,
w\}$ is a subset of the $d$-simplex $C$, it follows that $C_1$ is
a $(d - 1)$-simplex.

From the known structure of ${\rm lk}_{M_1}(A)$, we see that $D
:= A \cup (B \setminus \{y\}) \cup (R \setminus \{z^{\,\prime}\})
\cup (S \setminus \{w\})$ is also a $d$-simplex and hence $D_1 :=
V(M_1) \setminus \{x, y, z^{\,\prime}, w\} \subseteq D$ is a $(d -
1)$-simplex. One of the facet containing $D_1$ is $D$. Let
$D^{\,\prime}$ be the other facet containing $D_1$. If $x \in
D^{\,\prime}$ then $A \subseteq D^{\,\prime}$, but then from the
known link of $A$, we must have $D^{\,\prime} = D$, a
contradiction. So, $x \not\in D^{\,\prime}$.

If $z^{\,\prime} \in D^{\,\prime}$ then $C_1 \subseteq
D^{\,\prime}$. Since $C_1$ is also contained in the facets $C_1
\cup \{w\} = C$ and $C_1\cup \{x\}$, it follows that the
$(d-1)$-face $C_1$ is contained in three facets, a contradiction.
Thus, $z^{\,\prime} \not\in D^{\,\prime}$. If $y \in
D^{\,\prime}$ then $B \subseteq D^{\,\prime}$, which is not
possible by Claim 1. So, $y \not\in D^{\,\prime}$. Since $x, y,
z^{\,\prime}$ are three distinct vertices outside the facet
$D^{\,\prime}$ of $M_1$, and since $M_1$ is of excess 2, it
follows that $D^{\,\prime} = V(M_1) \setminus \{x, y,
z^{\,\prime}\} = C \cup \{z\} \setminus \{z^{\,\prime}\}$. This
proves Claim 2.

Now, in view of Claim 2, Lemma \ref{l8.1} $(i)$ and $(iii)$ imply
that there is a Gale diagram $G$ of $M_1$ such that all the
vertices in $R$ are represented by a single point $\rho$ and all
the vertices in $S$ are represented by a single point $\sigma$.
Let $z \in R$ and $w\in S$. Then, from the known link of $A$ in
$M_1$, one see that for any $y \in B$, $V(M_1)\setminus\{y, z,
w\}$ is a facet. Hence, by Lemma \ref{l4.2} $(a)$, $z \not\sim w$.
Therefore, by Lemma \ref{l8.1}, $\rho \neq \sigma$. Also, since
$\{y, z, w\}$ is a co-facet, it follows that $\{z, w\}$ is not a
co-facet. Hence $\sigma \neq \rho$. Thus, $\sigma \not\in \{\rho,
- \rho\}$. Again, since $\{y, z, w\}$ is a co-facet, it follows
that the point of $G$ representing $y \in B$ is in the shorter
open arc $(- \rho, - \sigma)$. Since $y \in B$ is arbitrary, this
shows that all the vertices in $B$ are represented by points in
the shorter open arc $(- \rho, - \sigma)$.

\smallskip

\noindent {\em Claim} 3. If $\gamma$ is a point in $G$
corresponding to a vertex $x \in A$, then $\gamma$ does not
belong to the shorter closed interval $[-\rho, -\sigma]$ of $G$.

\smallskip

If $\gamma = -\rho$ then for any $z \in R$, $\{z, x\}$ is a
co-facet in $M_1$. In that case, $V \setminus \{z, x\}$ is a facet
of $M_1$ containing $B$. This is not possible by Claim 1.
Therefore, $\gamma \neq - \rho$. Similarly, $\gamma \neq -
\sigma$. If $\gamma \in (-\rho, -\sigma)$, then for any $z\in R$
and $w\in S$, $\{x, z, w\}$ is a co-facet in $M_1$. Again this is
not possible by Claim 1. This proves this claim.

Let $\gamma \in (-\rho, \sigma)$ be the point of $G$ corresponding
to a vertex $x$ in $A$. If possible let there be a point $\mu$ of
$G$ in $(-\sigma, -\gamma]$. Let $y$ be a vertex corresponding to
this point. Clearly, $y\in A$. If $\mu = -\gamma$ then $\{y, x\}$
is a co-facet of $M_1$, which is not possible by Claim 1. If $\mu
\in (-\sigma, -\gamma)$ then for any $w\in S$, $\{x, y, w\}$ is a
co-facet of $M_1$. Again, this is not possible by Claim 1. Thus,
a point $\gamma$ of $G$ is in $(- \rho, \sigma)$ implies there is
no point of $G$ in $[- \sigma, - \gamma]$. Let $H_1$ be the Gale
diagram obtained from $G$ by shifting (all the occurrences)
$\gamma$ to $\sigma$. Then, by Lemma \ref{l8.1} $(iii)$, $H_1$ is
isomorphic to $G$ and hence is a Gale diagram of $M_1$. Continuing
this process (finitely many times) we get a Gale diagram $H_2$ of
$M_1$ in which points corresponding to the vertices in $B \cup R
\cup S$ are same as those in $G$ and the points corresponding to
the vertices of $A$ are in $[\sigma, \rho] \cup (\rho, -\sigma)$.
Since there is no points of $H_2$ in $(-\rho, \sigma)$, by a
similar process we get a Gale diagram $H_3$ of $M_1$ in which
points corresponding to the vertices $B \cup R \cup S$ are same
as those in $G$ and the points corresponding to the vertices of
$A$ are in $[\sigma, \rho]$.

If the vertices of $A$ correspond to points in $(\rho, \sigma)$
then $G_1 = H_3$ is a Gale diagram as required. Otherwise, let
$\rho$ correspond to some vertices of $A$ as well.

Since there is no points of $H_3$ at $- \rho$, by Lemma
\ref{l8.1} $(iii)$, $H_3$ is isomorphic to a Gale diagram $H_4$
in which there is one more point, say $\nu$, than $H_3$ near
$\rho$ in $(\rho, \sigma)$. Then $H_4$ is a Gale diagram of $M_1$
in which the point $\nu$ corresponds to those vertices of $A$
which correspond to $\rho$ in $H_3$. So, $\rho$ (in $H_4$) does
not correspond to any vertex of $A$. By doing similar process, we
get a Gale diagram $H_5$ in which $\sigma$ does not correspond to
any vertex of $A$. Now, $G_1 = H_5$ is a Gale diagram as required.
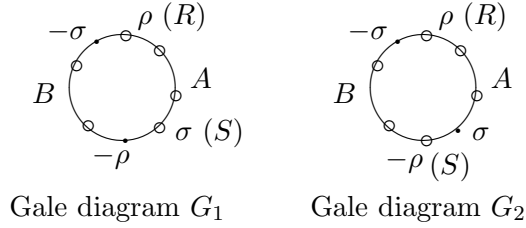
\hfill $\Box$

\setlength{\unitlength}{1mm}

\begin{picture}(80,35)(-40,-13)


\put(15,8){\circle{15}}

\put(14.5,14){$\circ$} \put(14.8,0.7){\bf .}

\put(19,1.7){$\circ$} \put(11,13.8){\bf .}

\put(9.5,2){$\circ$} \put(8,10){$\circ$}

\put(21.2,6){$\circ$} \put(19,12){$\circ$}

\put(17,16.5){\mbox{$\rho$ $(R)$}} \put(11,-2){\mbox{$-\rho$}}
\put(5,14.5){\mbox{$-\sigma$}} \put(22,1){\mbox{$\sigma$ $(S)$}}
\put(3,6){\mbox{$B$}} \put(24,8){\mbox{$A$}}

\put(0,-9){\mbox{Gale diagram $G_1$}}

\put(55,8){\circle{15}}

\put(54.5,14){$\circ$} \put(54.5,0){$\circ$}

\put(59,2){\bf .} \put(51,13.8){\bf .}

\put(49.5,2){$\circ$} \put(48,10){$\circ$}

\put(61.2,6){$\circ$} \put(59,12){$\circ$}

\put(57,16.5){\mbox{$\rho$ $(R)$}} \put(50,-2.5){\mbox{$-\rho$}}
\put(55.8,-3.5){\mbox{$(S)$}} \put(45,14.5){\mbox{$-\sigma$}}
\put(61.5,1){\mbox{$\sigma$}} \put(43,6){\mbox{$B$}}
\put(64,8){\mbox{$A$}}

\put(40,-9){\mbox{Gale diagram $G_2$}}

\end{picture}

\noindent {\em Proof of Lemma} \ref{l8.3}. Let $K_1 = S^{\,a}_{a
+ 2}(A) \ast S^{\,b}_{b + 2}(B) \ast S^{\,r}_{r + 2}(R) \ast
S^{\,s}_{s + 2}(S)$. Let $u\not\in V(M)$. For $2\leq i\leq 4$,
consider the pure $d$-dimensional simplicial complex $K_i$ with
facet-set $F(K_i)$ given by:
\begin{eqnarray*}
F(K_2) & = & \{A \cup (B \setminus \{y\}) \cup (R\setminus \{z\})
\cup (S\setminus\{w\}) : y \in B, z \in R, w \in S\}, \\
F(K_3) & = & \{\{u\} \cup \alpha : \alpha \mbox{ is a facet of }
K_1\}, \\
F(K_4) & = & \{\{u\} \cup \beta : \beta \mbox{ is a facet of }
B(K_1)\}.
\end{eqnarray*}
Observe that both $B(K_2)$ and $K_4$ are triangulations of a
$d$-ball with boundary $B(K_1)$.

Clearly, $K_3$ is obtained from $K_2$ by starring the vertex $u$
in the face $A$. So, $K_3$ is a subdivision of $K_2$ and $K_4$ is
a subdivision of $K_3$. So, $K_4$ is a subdivision of $K_2$.
Therefore, $K_4$ is combinatorially equivalent to $B(K_2)$.

 Since $B(K_2) \subseteq B(M_1)$,
$K_4 \subseteq B(M)$ and $B(M_1)$ is obtained from $B(M)$ by
replacing $K_4$ by $B(K_2)$, $B(M_1)$ is equivalent to $B(M)$.
This proves $(i)$.

$(ii)$ follows from $(i)$.

To prove $(iii)$, let $G_1$ be as in Lemma \ref{l8.4} and let $V =
V(M) = V(M_1)$.

Consider the polytopal sphere $M_2$ with vertex-set $V$ and Gale
diagram $G_2$ obtained from $G_1$ by replacing $\sigma$ by
$\rho$. (This is a Gale diagram since $\#(A) \geq 2$.) The
correspondence between the vertex-set $V$ of $M_2$ and the
multi-set $G_2$ of points is as follows. The vertices in $A$ and
$B$ are represented by the same points of $G_2$ as in $G_1$. The
vertices in $R$ are represented by $\rho$ and the vertices in $S$
are represented by $- \rho$.

Let $\alpha$ be a facet of $M$. If $\tau \subseteq \alpha$ then
$\alpha$ is of the form $\tau \cup (R \setminus \{z\}) \cup
(S\setminus \{w\})$ for some $z \in R$ and $w \in S$. Since the
points corresponding to $z$ and $w$ in $G_2$ are $\rho$ and $-
\rho$ respectively, $\sigma = V \setminus \{z, w\}$ is a facet of
$M_2$.

If $\tau \not \subseteq \alpha$ then $\alpha$ is a facet of $M_1$.
Since $\tau$ is the smallest face of $M$ containing $A$ and $\tau
\not \subseteq \alpha\in M$, $A \not \subseteq \alpha$. Similarly,
$B \not \subseteq \alpha$. Thus, if $\alpha^{\,\prime} := V
\setminus \alpha$ then $\alpha^{\,\prime} \cap A \neq \emptyset$
and $\alpha^{\,\prime} \cap B \neq \emptyset$. Let $x\in A\cap
\alpha^{\,\prime}$ and $y\in B \cap \alpha^{\,\prime}$. Let the
points corresponding to $x$ and $y$ in $G_1$ be $\gamma\in (\rho,
\sigma)$ and $\mu \in (-\rho, -\sigma)$ respectively.

If $\alpha^{\,\prime}$ has two vertices then $\alpha^{\,\prime} =
\{x, y\}$.  Since $\alpha\in M_1$, $\mu = -\gamma$ (in $G_1$ and
in $G_2$) and hence $\{x, y\}$ is a co-facet in $M_2$. Therefore,
$\alpha\in M_2$. If $\alpha^{\,\prime}$ has three vertices then
$\alpha^{\,\prime} = \{x, y, u\}$. Let the point corresponding to
$u$ in $G_1$ be $\nu$. Then the interior of the convex hull of
$\{\gamma, \mu, \nu\}$ contains the origin $0$. If $\nu \neq
\sigma$ then $\gamma, \mu, \nu \in G_2$ and hence $\{x, y, u\}$
is a co-facet in $M_2$. Then $\alpha$ is a facet of $M_2$. If
$\nu = \sigma$ then, since $\alpha$ is a facet of $M_1$, the
interior of the convex hull of $\{\gamma, \mu, \sigma\}$ contains
$0$. In this case, since $\gamma\in (\rho, \sigma)$ and $\mu \in
(-\rho, -\sigma)$, it is easy to see that the interior of the
convex hull of $\{\gamma, \mu, - \rho\}$ also contains $0$. Since
$u$ corresponds to $- \rho$ in $G_2$, $\{x, y, u\}$ is a co-facet
of $M_2$. So, $\alpha$ is a facet of $M_2$. Therefore, $M$ is a
sub-pseudomanifold of $M_2$ and hence, $M = M_2$ of the same
dimension. This implies that $M$ is a polytopal sphere. This
completes the proof. \hfill $\Box$

\begin{theo}  \label{t8.1}
Every cellular sphere of excess $2$ is a polytopal sphere.
\end{theo}

\noindent {\em Proof.} Let $M$ be a cellular $d$-sphere with
$d+4$ vertices.

\noindent {\em Case} 1\,: Assume that $M$ has no facet containing
$d+3$ vertices. We may also assume that $M$ has no face whose
boundary is the join of two standard spheres. (If we have the
result in this case then Lemma \ref{l8.3} enables us to handle
the general case by a simple induction on the number of faces
whose boundary is such a join.) Then, by Lemma \ref{l8.2} $(ii)$,
$M$ is a combinatorial sphere. Since all combinatorial spheres of
excess at most 2 are polytopal spheres (cf. \cite{bd2}), it
follows that $M$ is a polytopal sphere.

\noindent {\em Case} 2\,: Now, assume that $M$ has a facet $\tau$
containing $d+3$ vertices. Then $M = (\partial_M \tau) \otimes
S^{\,-1}_1$. Let $c$ be the largest integer such that $M = N
\otimes S^{\,c}_{c + 2}$ for some cellular pseudomanifold $N$.
Then, by Lemma \ref{l2.5}, the dimension of $N$ is $d- c -2$,
$e(N) = 2$ and $N$ has no facet containing $d-c+1$ vertices (else
$N = N_0 \otimes S^{\,-1}_1$ for some cellular pseudomanifold
$N_0$ and hence, by Lemma \ref{l2.8} $(a)$, $M = N_0 \otimes
S^{\,c + 1}_{c +3}$, contradicting the choice of $c$). Then, by
Lemma \ref{l2.8} $(d)$, $N$ is a cellular sphere. So, $N$ is a
cellular $(d- c -2)$-sphere of excess 2 and has no facet
containing $d-c+1$ vertices. Therefore, by Case 1, $N$ is a
polytopal sphere and hence by Lemma \ref{l2.8} $(e)$, $M$ is a
polytopal sphere. This completes the proof. \hfill $\Box$

\end{document}